\documentclass[11 pt]{amsart}
\usepackage{import}

\usepackage{amsaddr}
\usepackage{amssymb}
\usepackage{asymptote}
\usepackage{chngcntr}
\usepackage{enumitem}
\counterwithin{figure}{section}
\counterwithin{table}{section}
\numberwithin{equation}{section}
\usepackage{float}
\usepackage[margin=1.5in]{geometry}
\usepackage{graphicx}
\usepackage{hyperref}
\usepackage[parfill]{parskip}
\usepackage{tablefootnote}
\usepackage{tikz-cd}
\usepackage{todonotes}
\usepackage[enableskew,vcentermath]{youngtab}

\makeatletter
\newtheorem*{rep@theorem}{\rep@title}
\newcommand{\newreptheorem}[2]{%
\newenvironment{rep#1}[1]{%
 \def\rep@title{#2 \ref{##1}}%
 \begin{rep@theorem}}%
 {\end{rep@theorem}}}
\makeatother

\theoremstyle{plain}
\newtheorem{thm}{Theorem}[section]
\newtheorem{lem}[thm]{Lemma}
\newtheorem{cor}[thm]{Corollary}
\newtheorem{prop}[thm]{Proposition}
\newreptheorem{thm}{Theorem}

\theoremstyle{definition}
\newtheorem{defn}[thm]{Definition}
\newtheorem{ex}[thm]{Example}
\newtheorem{prob}[thm]{Problem}

\theoremstyle{remark}
\newtheorem{rem}[thm]{Remark}

\newcommand{\vocab}[1]{\textit{#1}}

\newcommand{\lt}{\left}
\newcommand{\rt}{\right}

\newcommand\blank{$ $}
\newcommand\bull{\bullet}
\newcommand\rcirc{\color{red}{\circ}}
\newcommand\btimes{\color{blue}{\times}}
\newcommand\ten{10}
\newcommand\eleven{11}
\newcommand\twelve{12}
\newcommand\thirteen{13}
\newcommand\fourteen{14}
\newcommand\fifteen{15}
\newcommand\sixteen{16}
\newcommand\seventeen{17}
\newcommand\eighteen{18}
\newcommand\nineteen{19}
\newcommand\twenty{20}
\newcommand\twone{21}
\newcommand\twtwo{22}
\newcommand\twthree{23}
\newcommand\twfour{24}
\newcommand\twfive{25}
\newcommand\twsix{26}
\newcommand\twseven{27}
\newcommand\tweight{28}
\newcommand\twnine{29}
\newcommand\rbsix{{\color{red}{\mathbf {6}}}}
\newcommand\rbnine{\color{red}{\mathbf {9}}}
\newcommand\rbten{\color{red}{\mathbf {10}}}
\newcommand\rbeleven{\color{red}{\mathbf {11}}}
\newcommand\rbtwelve{\color{red}{\mathbf {12}}}
\newcommand\rbthirteen{\color{red}{\mathbf {13}}}
\newcommand\rbfourteen{\color{red}{\mathbf {14}}}
\newcommand\rbfifteen{\color{red}{\mathbf {15}}}

\newcommand\bone{\mathbf {1}}
\newcommand\btwo{\mathbf {2}}
\newcommand\bfour{\mathbf {4}}
\newcommand\bfive{\mathbf {5}}
\newcommand\bsix{\mathbf {6}}
\newcommand\beight{\mathbf {8}}
\newcommand\bnine{\mathbf {9}}
\newcommand\bten{\mathbf {10}}
\newcommand\beleven{\mathbf {11}}
\newcommand\btwelve{\mathbf {12}}
\newcommand\bthirteen{\mathbf {13}}
\newcommand\bfourteen{\mathbf {14}}
\newcommand\bfifteen{\mathbf {15}}

\newcommand{\padTableText}[1]{ {\begin{tabular}{c} #1 \end{tabular}} }
\newcommand{\younganddescent}[2]{
    \text{
        \begin{tabular}{l}
            $\young(#1)$ \\[10px]
            #2
        \end{tabular}
    }
}
\newcommand{\padYD}[2]{
    \begin{tabular}{c}
        \\[-8px]
        #1 \\ [#2 px]
    \end{tabular}
}

\newcommand{\skewshape}{{\lambda/\mu}}
\newcommand{\ZZ}{{\mathbb{Z}}}
\newcommand{\cS}{{\mathcal{S}}}
\newcommand{\functo}{{~\rightarrow~}}
\newcommand{\phimapsto}{{\overset{\phi}{\mapsto}}}
\newcommand{\jdt}{\widetilde{\rm{jdt}}}
\newcommand{\RotSE}{\Rot_{\SE}}
\newcommand{\RotNW}{\Rot_{\NW}}
\newcommand{\RcSE}{\Rc_{\SE}}
\newcommand{\RcNW}{\Rc_{\NW}}
\newcommand{\RpSE}{\Rp_{\SE}}
\newcommand{\RpNW}{\Rp_{\NW}}
\newcommand{\RpSEc}{\RpSE^{\cen}}
\newcommand{\RpSEsw}{\RpSE^{\SW}}
\newcommand{\RpSEne}{\RpSE^{\NE}}

\DeclareMathOperator{\Jdt}{jdt}
\DeclareMathOperator{\Des}{Des}
\DeclareMathOperator{\cDes}{cDes}
\DeclareMathOperator{\SYT}{SYT}
\DeclareMathOperator{\pos}{pos}
\DeclareMathOperator{\pro}{pro}
\DeclareMathOperator{\dem}{dem}
\DeclareMathOperator{\tra}{t}
\DeclareMathOperator{\rev}{rev}
\DeclareMathOperator{\Rot}{rot}
\DeclareMathOperator{\Rc}{Rc}
\DeclareMathOperator{\Rp}{Rp}
\DeclareMathOperator{\SE}{SE}
\DeclareMathOperator{\NW}{NW}
\DeclareMathOperator{\cen}{C}
\DeclareMathOperator{\SW}{SW}
\DeclareMathOperator{\NE}{NE}

\begin{document}

\title[Cyclic Descents]{Cyclic Descents for General Skew Tableaux}

\author{Brice Huang}
\address{Department of Mathematics, MIT, Cambridge, MA, USA.}
\email{bmhuang@mit.edu}


\date{August 14, 2018}

\begin{abstract}
A \vocab{cyclic descent} function on standard Young tableaux of size $n$
is a function that restricts to the usual descent function when $n$ is omitted,
such that the number of standard Young tableaux of given shape
with cyclic descent set $D\subset[n]$
is invariant under any modulo $n$ shift of $D$.
The notion of cyclic descent was first studied for rectangles by Rhoades,
and then generalized to certain families of skew shapes
by Adin, Elizalde, and Roichman.
Adin, Reiner, and Roichman proved that a skew shape has a cyclic descent
map if and only if it is not a connected ribbon.
Unfortunately, their proof is nonconstructive;
until now, explicit cyclic descent maps are known only for
small families of shapes.

In this paper, we construct an explicit cyclic descent map for
all shapes where this is possible.
We thus provide a constructive proof of Adin, Reiner, and Roichman's
result.
Our construction of a cyclic descent map generalizes many of the constructions
in the literature.
\end{abstract}

\maketitle

\section{Introduction}\label{sec:intro}

Let $\cS_n$ denote the symmetric group on $[n]=\{1,\ldots,n\}$.
For a permutation $\pi\in \cS_n$, the \vocab{descent set} of $\pi$ is defined by
\begin{equation*}
    \Des(\pi) = \{i\in [n-1]: \pi(i) > \pi(i+1)\}.
\end{equation*}
In 1998, Cellini \cite{Cel} introduced a notion of \vocab{cyclic descent}.
Cyclic descents were further studied in \cite{Pet}, and are defined as follows.
For $\pi \in \cS_n$, let
\begin{equation*}
    \cDes(\pi) =
    \begin{cases}
        \Des(\pi) \cup \{n\} & \pi(n) > \pi(1), \\
        \Des(\pi) & \pi(n) < \pi(1).
    \end{cases}
\end{equation*}

The function $\cDes$ has the property that the multiset
\begin{equation*}
    \{\{\cDes(\pi) | \pi \in \cS_n\}\}
\end{equation*}
is invariant under rotation of indices modulo $n$.
Equivalently, the number of permutations with cyclic descent $D\subseteq [n]$ is
invariant under any modulo $n$ rotation of $D$.

Moreover, let $\phi: \cS_n \functo \cS_n$ be the cyclic rotation function,
such that
\begin{equation*}
    (\phi \pi) (i) = \pi (i-1)
\end{equation*}
for each $\pi \in \cS_n$, where indices are taken modulo $n$.
Then $\phi$ has the property that for any $\pi \in \cS_n$,
\begin{equation*}
    \cDes(\pi)+1 = \cDes(\phi \pi),
\end{equation*}
where $\cDes(\pi)+1$ denotes the set obtained from $\cDes(\pi)$
by incrementing each element by $1$ modulo $n$.

In fact, the existence of a bijection $\phi$ with this property
implies that the multiset $\{\{\cDes(\pi) | \pi \in \cS_n\}\}$
is rotation-invariant modulo $n$.

Standard Young tableaux also have a well-studied notion of descent.
Throughout this paper, let $\skewshape$ denote a skew shape with
$|\skewshape| = n$, where $\mu$ is a Young diagram contained in $\lambda$;
we draw Young diagrams with origin in the northwest corner,
so that in each Young diagram,
the rows and columns are justified on the west and north, respectively.

The \vocab{standard Young tableaux} of shape $\skewshape$ are the labellings of
the cells of $\skewshape$ with a permutation of $[n]$,
such that row entries are increasing from north to south,
and column entries are increasing from west to east.
Let $\SYT(\skewshape)$ denote the set of standard Young tableaux of shape
$\skewshape$.
For $T\in \SYT(\skewshape)$ the \vocab{descent set} of $T$ is given by
\begin{equation*}
    \Des(T) = \{
        i\in [n-1]: \text{the row of $i+1$ is strictly south of the row of $i$}
    \}.
\end{equation*}

\begin{ex}
    The standard Young tableau
    \begin{equation*}
        T = \young(::35,147,26)
    \end{equation*}
    of shape $\skewshape = (4,3,2)/(2)$ has descent set $\Des(T) = \{1,3,5\}$.
\end{ex}

A natural question is whether an analogous notion of cyclic descent
exists for standard Young tableaux.
This question was first introduced by Rhoades \cite{Rho} in 2010 for rectangular
tableaux, in the context of the \vocab{cyclic sieving phenomenon}
of Reiner, Stanton, and White \cite{RSW}.
For more on cyclic sieving, see \cite{AKLM, FK, Pec, PPR}.

The following formulation of cyclic descent is due to Adin, Reiner,
and Roichman in 2017.

\begin{defn}\cite[Definition 2.1]{ARR}\label{defn:cyc-descent-map}
    A \vocab{cyclic descent map} for $\skewshape$ is a pair $(\cDes, \phi)$,
    where $\cDes: \SYT(\skewshape) \functo 2^{[n]}$ is a function and
    $\phi: \SYT(\skewshape) \functo \SYT(\skewshape)$ is a bijection,
    that satisfies the following properties for all $T\in \SYT(\skewshape)$.
    \begin{itemize}
        \item (Extension) $\cDes(T)\cap [n-1] = \Des(T)$.
        \item (Equivariance) $\cDes(\phi T) = \cDes(T)+1$, with indices taken
        modulo $n$.
        \item (Non-Escher) $\emptyset \subsetneq \cDes(T) \subsetneq [n]$.
    \end{itemize}
\end{defn}
Note that the term ``cyclic descent map" exclusively
refers to the pair $(\cDes, \phi)$, not the cyclic descent function $\cDes$.
Moreover, note that $\cDes$ is uniquely determined by specifying, for each
$T\in \SYT(\skewshape)$, whether $n\in \cDes(T)$.

\begin{ex}
    The shape
    \begin{equation*}
        \skewshape = \young(\blank\blank\blank,\blank\blank\blank)
    \end{equation*}
    admits the following cyclic descent map $(\cDes, \phi)$.
    For each $T\in \SYT(\skewshape)$, the set $\cDes(T)$ is shown below it.
    \begin{equation*}
        \younganddescent{134,256}{$\{1,4\}$}
        \younganddescent{125,346}{$\{2,5\}$}
        \younganddescent{123,456}{$\{3,\rbsix\}$}
        \younganddescent{135,246}{$\{1,3,5\}$}
        \younganddescent{124,356}{$\{2,4,\rbsix\}$}
    \end{equation*}
    In constructing this map, we chose to assign $6$
    (shown in \textcolor{red}{\bf red}) to the cyclic descent sets
    of the third and fifth tableaux, as shown.
    For this particular $\skewshape$, this assignment is the only valid choice;
    in general there may be multiple valid choices,
    corresponding to multiple valid functions $\cDes$.
    The bijection $\phi$ cycles the first three tableaux and last two tableaux.
    Thus, $\skewshape$ has a cyclic descent map.
\end{ex}

This definition motivates the following problem.
\begin{prob}\cite{AER}
    For which skew shapes $\skewshape$ does a cyclic descent map exist?
\end{prob}
This problem was fully, albeit nonconstructively, solved by Adin, Reiner, and Roichman
in \cite{ARR} using nonnegativity properties of Postnikov's toric Schur polynomials \cite{Pos}.

A \vocab{connected ribbon} is a connected skew shape
with no $2 \times 2$ square.
\begin{thm}\cite[Theorem 1.1]{ARR}\label{thm:arr}
    Let $\skewshape$ be a skew shape with $n$ cells.
    A cyclic descent map exists for $\skewshape$ if and only if $\skewshape$ is not a connected ribbon.
    Moreover, for all $J\subseteq [n]$, all cyclic descent maps share the same fiber sizes $|\cDes^{-1}(J)|$.
\end{thm}
However, the proof does not explicitly construct the cyclic descent map.
It would be desirable, therefore, to prove Theorem~\ref{thm:arr} combinatorially.
Adin, Reiner, and Roichman \cite{ARR} give a combinatorial proof that
connected ribbons do not have cyclic descent maps.
They also combinatorially prove that the fiber sizes $|\cDes^{-1}(J)|$ are fixed
and give an explicit formula for the fiber sizes.

So, it remains to solve the following problem, posed by Adin, Elizalde, and Roichman in 2018.
This problem reformulates \cite[Problems 7.1, 7.2]{AER}.

\begin{prob}\cite{AER}\label{prob:main}
    For $\skewshape$ not a connected ribbon, explicitly construct a cyclic descent map $(\cDes, \phi)$.
\end{prob}

It is worth noting that solving this problem provides a combinatorial procedure
for computing certain Gromov-Witten invariants that appear as structure constants
in the quantum cohomology of Grassmannians \cite{ARR}.
These Gromov-Witten invariants have many combinatorial and algebraic interpretations,
summarized in \cite[Section 3.4]{ARR}.

There are many partial results on this problem in the literature,
the first of which is due to Rhoades \cite{Rho}.
Let $\pro$ and $\dem$ be the (Sch\"utzenberger) promotion and demotion operators, respectively;
these will be formally defined in Section~\ref{sec:prelim}.\footnote{
    The literature contains two inconsistent definitions of promotion.
    \cite{AER,StaEC1,StaPE} define promotion to be the operation that we term demotion.
    We follow the definition from \cite{ER, Rho}.
}

\begin{thm}\label{thm:rhoades}\cite[Lemma 3.3]{Rho}
    Let $\skewshape$ be a rectangular Young diagram with length and width both larger than 1.
    Let $\phi = \pro$, and let $n\in \cDes(T)$ if and only if $n-1 \in \Des(\dem T)$.
    Then $(\cDes, \phi)$ is a cyclic descent map.
\end{thm}

\begin{rem}
    By definition of cyclic descent map,
    and because $\phi = \pro$ and $\phi^{-1} = \dem$,
    the condition $n-1 \in \Des(\dem T)$ is equivalent to $1\in \Des(\pro T)$.
    The latter condition is more analogous to our main result;
    we state Theorem~\ref{thm:rhoades} using the former condition
    to preserve the original statement.
\end{rem}

Analogous results exist in the literature for other families
of skew diagrams $\skewshape$.  A summary of these results
can be found in Table~\ref{table:literature-summary} below.

\begin{table}[H]
    \centering
    \begin{tabular}{|c|c|c|}
        \hline
        Shape & Example & Reference \\
        \hline
        Rectangle with dimensions $>1$ &
        \padYD{$\young(\blank\blank\blank,\blank\blank\blank)$}{10} &
        \cite[Lemma 3.3]{Rho} \\
        \hline
        \padTableText{Young diagram and \\ disconnected northeast cell} &
        \padYD{$\young(:::\blank,\blank\blank\blank,\blank\blank)$}{17} &
        \cite[Proposition 5.3]{ER}  \\
        \hline
        Hook plus internal cell &
        \padYD{$\young(\blank\blank\blank\blank,\blank\blank,\blank)$}{17} &
        \cite[Theorem 1.11]{AER} \\
        \hline
        Strip\tablefootnote{
            A disconnected union of single-row or single-column rectangles.
        } &
        \padYD{$\young(:::\blank\blank\blank,::\blank,::\blank,\blank\blank)$}{24} &
        \cite[Proposition 3.3]{AER} \\
        \hline
        \padTableText{2-row straight shape,\\ not ribbon} &
        \padYD{$\young(\blank\blank\blank\blank,\blank\blank\blank)$}{10} &
        \cite[Theorem 1.14]{AER} \\
        \hline
        \padTableText{2-row skew shape,\\ not ribbon} &
        \padYD{$\young(:\blank\blank\blank\blank,\blank\blank\blank)$}{10} &
        \cite[Theorems 6.3, 6.11]{AER}\tablefootnote{
            Unlike the above results, these two theorems only give maps $\cDes$
            such that $\{\{\cDes(T) | T\in \SYT(\skewshape)\}\}$
            is invariant under modulo $n$ rotation.
            They do not give explicit $\phi$,
            and explicit $\phi$ corresponding to these $\cDes$ are not known.
        }\\
        \hline
    \end{tabular}
    \caption{Summary of literature on Problem~\ref{prob:main}.}
    \label{table:literature-summary}
\end{table}

In this paper, we answer Problem~\ref{prob:main} completely.
We construct a general cyclic descent map for all $\skewshape$ other than
connected ribbons.
This construction requires two bijective operations on $\SYT(\skewshape)$,
which we call \vocab{southeast rotation} and \vocab{northwest rotation},
and which we denote $\RotSE$ and $\RotNW$.
These operations are defined in Section~\ref{sec:rotation}
and are new to the literature.

Our main result is as follows.
\begin{thm}\label{thm:main}
Suppose $\skewshape$ is a skew shape that is not a connected ribbon.
Let
\begin{equation*}
    \phi = \RotNW^{-1} \circ \pro \circ \RotSE
\end{equation*}
and let $n\in \cDes(T)$ if and only if $1\in \Des(\phi T)$.
Then $(\cDes, \phi)$ is a cyclic descent map.
\end{thm}
Let us provide some motivation for the rotation operators $\RotSE$ and $\RotNW$.
The \vocab{promotion path} of a standard Young tableau $T$ be the set of cells
traversed by $1$ when $T$ is promoted.
One proof of Theorem~\ref{thm:rhoades} uses the following property of the promotion paths.
For a standard Young tableau $T$ of any shape, let $P_1$ and $P_2$ be the promotion paths of $T$ and $\pro T$,
and define the set of cells
\begin{equation*}
    \Gamma = \{Z \in P_1 | \text{western neighbor of $Z$ is also in $P_1$}\}.
\end{equation*}
The property is that $P_2$ does not cross $P_1$ from south to north,
in the sense that if $Z\in \Gamma \cap P_2$,
then the southern neighbor of $Z$ is not in $P_2$.

Suppose now $T$ has rectangular shape.
The proof uses the non-crossing property to show that $n-1\in \Des(T)$ implies $1\in \Des(\pro^2 T)$.
If $n-1\in \Des(T)$, the first step of $P_1$ is northward, and
the non-crossing property ensures that $P_2$ is trapped south of $P_1$.
This forces the last step of $P_2$ to be northward, which implies that $1\in \Des(\pro^2 T)$.
We symmetrically get that $n-1\not\in \Des(T)$ implies $1\not\in \Des(\pro^2 T)$,
so $n-1\in \Des(T)$ if and only if $1\in \Des(\pro^2 T)$.
This implies equivariance, which is the core of the proof of Theorem~\ref{thm:rhoades}.\footnote{
This strategy is equivalent to
the strategy used by Rhoades in \cite{Rho}, though Rhoades's argument,
which uses demotion instead of promotion, is formulated slightly differently.}

When $T$ does not have rectangular shape,
this strategy fails because there are now multiple possible
source and destination cells for $P_1$ and $P_2$.
The fact that $n-1\in \Des(T)$ no longer implies $P_2$ is trapped south of $P_1$,
because $P_2$ may start at a more northern cell than $P_1$;
moreover, even when $P_2$ is trapped south of $P_1$, the last step of $P_2$ is not necessarily
northward, because $P_2$ may end at a different cell than $P_1$.

This is the problem the rotation operators solve.
The rotation operators force the source and destination of $P_2$ to be,
respectively, southwest and northeast of the source and destination of $P_1$;
due to this relative positioning, the non-crossing property once again
traps $P_2$ south of $P_1$, and the above strategy succeeds.
This argument is carried out in the proof of Proposition~\ref{prop:double-promotion}.
It emulates the proof of Theorem~\ref{thm:rhoades}, though in substantially generalized form.

As we will see in Section~\ref{sec:special-cases},
Theorem~\ref{thm:main} generalizes most of the constructions in
Table~\ref{table:literature-summary}.

The proof of Theorem~\ref{thm:main} also explains the need for the hypothesis
that $\skewshape$ is not a connected ribbon.
This hypothesis is necessary for the crucial Lemma~\ref{lem:rot-no-interference},
which states that the two rotation operations
$\RotSE$ and $\RotNW$, in an appropriate sense, do not interfere.
That this non-interference property does not hold
when $\skewshape$ is a connected ribbon motivates the need for this hypothesis.

The rest of this paper is structured as follows.
In Section~\ref{sec:prelim}, we review the definitions of promotion and
demotion, state their relevant properties, and introduce two symmetry operations
on standard Young tableaux.
In Section~\ref{sec:rotation}, we define the rotation operations and prove their
key properties.
Section~\ref{sec:discussion} presents some examples of Theorem~\ref{thm:main}
and discusses some properties of the construction $(\cDes, \phi)$.
Section~\ref{sec:main-proof} is devoted to the proof of Theorem~\ref{thm:main}.
Section~\ref{sec:special-cases}
discusses the cyclic descent maps known in the literature
in relation to Theorem~\ref{thm:main};
most of these constructions are special cases of Theorem~\ref{thm:main}.
Finally, Section~\ref{sec:conclusion} concludes with some open problems.

\subsection*{Acknowledgments}

This research was funded by NSF/DMS grant 1650947 and NSA grant H98230-18-1-0010
as part of the 2018 Duluth Research Experience for Undergraduates (REU) program.
The author thanks Joe Gallian for supervising the research and suggesting the problem.
The author also thanks Joe Gallian, Ashwin Sah, and Noah Kravitz for close
readings of drafts of this paper,
and Aaron Berger and Colin Defant for useful discussions.

\section{Preliminaries}\label{sec:prelim}

\subsection{Promotion and Demotion}

We first define the Sch\"utzenberger promotion and demotion operators,
which will appear throughout this paper.
\begin{defn}
    The action of the \vocab{promotion} operator $\pro$
    on $T\in \SYT(\skewshape)$ is as follows.
    Add $1 \pmod n$ to each entry of $T$;
    this turns the $n$ into a $1$ and increments the remaining entries.
    Repeatedly apply the following operation until the cell containing $1$
    has neither northern nor western neighbor:
    swap $1$ with the larger of its northern neighbor (if it exists)
    and its western neighbor (if it exists).
\end{defn}
It is not difficult to show that the resulting tableau is standard.

Analogously, the \vocab{demotion}
(sometimes known in the literature as \vocab{dual promotion})
operator is defined by subtracting $1 \pmod n$ from each entry of $T$,
and repeatedly swapping $n$ with the smaller of its southern and eastern neighbors.
It is clear that promotion and demotion are inverses.

Promotion and demotion can be equivalently defined using the
\vocab{jeu de taquin} ($\Jdt$) operation,
a formal definition of which can be found in \cite[p. 419-420]{StaEC2}.
To promote a tableau $T\in \SYT(\skewshape)$,
replace the cell containing $n$ with a $\Jdt$ hole,
perform a series of $\Jdt$ slides to move the hole to a northwestern corner,
increment all entries, and fill in the hole with a cell with entry $1$.
To demote a tableau $T\in \SYT(\skewshape)$,
replace the cell containing $1$ with a $\Jdt$ hole, decrement all entries,
move the hole to a southeastern corner with $\Jdt$ slides,
and fill the hole with a cell with entry $n$.

Jeu de taquin has the following property.
\begin{lem}\label{lem:jdt-preserve-des}\cite[p. 431]{StaEC2}
    A series of $\Jdt$ slides preserves the descent set of any standard Young tableau.
\end{lem}

The following property of promotion and demotion will be useful.
\begin{lem}\label{lem:pro-increments-des}
    Let $\skewshape$ be any skew shape.
    For $i\in \{1,\ldots,n-2\}$, $i\in \Des(T)$ if and only if $i+1 \in \Des(\pro T)$.
    Similarly, for $i\in \{2,\ldots,n-1\}$, $i\in \Des(T)$ if and only if $i-1\in \Des (\dem T)$.
\end{lem}
\begin{proof}
    Let us consider the steps in the procedure of the jeu de taquin definition of promotion.
    Deleting the cell containing $n$ does not modify any descents in $\{1,\ldots,n-2\}$.
    By Lemma~\ref{lem:jdt-preserve-des}, sliding the $\Jdt$ hole to the northwest by
    $\Jdt$ slides also preserves the descents in $\{1,\ldots,n-2\}$.
    Incrementing all the entries increments these descents,
    so for each descent $i\in \{1,\ldots,n-2\}$, $i+1$ is a descent of the incremented tableau.
    Filling the hole with a cell with entry $1$ does not affect these descents.
    This proves the lemma for promotion.

    The argument for demotion is analogous.
\end{proof}

\begin{ex}
    For $\lambda = (5,3,2,2,1)$ and $\mu = (2,2)$,
    below is a tableau $T\in \SYT(\skewshape)$ and its promotion.
    \begin{equation*}
        \young(::268,::3,14,59,7)
        \quad
        \overset{\pro}{\mapsto}
        \quad
        \young(::379,::4,15,26,8).
    \end{equation*}
    Note that $\Des(T) = \{2,3,4,6,8\}$ and $\Des(\pro T) = \{1,3,4,5,7\}$,
    so $i\in \Des(T)$ if and only if $i+1\in \Des(\pro T)$ for $i\in \{1,\ldots,7\}$.
\end{ex}

\subsection{Transpose and Reverse}

The \vocab{transpose} of a skew shape $\skewshape$,
denoted $(\skewshape)^{\tra}$,
is the skew shape obtained by reflecting $\skewshape$ over a
northwest-southeast line.
This operation interchanges rows with columns.
For $T\in \SYT(\skewshape)$, the transpose of $T$, denoted $T^{\tra}$,
is the reflection of $T$ over a northwest-southeast line.

The \vocab{reverse} of a skew shape $\skewshape$, denoted $(\skewshape)^{\rev}$,
is the shape obtained by rotating $\skewshape$ by $180^\circ$.
The reverse of a tableau $T\in \SYT(\skewshape)$, denoted $T^{\rev}$,
is the tableau obtained by rotating $T$ by $180^\circ$
and replacing each entry $i$ by $n+1-i$.

The \vocab{reverse transpose}
of a skew shape $\skewshape$ or a tableau $T\in \SYT(\skewshape)$
is the shape obtained from it by applying the reverse and transpose operations;
note that these operations commute, so the reverse transpose is well defined.

The following properties of these operators are clear.
Promotion and demotion commute with transposition,
and are conjugate with respect to reversal.
That is, for all $T\in \SYT(\skewshape)$,
\begin{eqnarray*}
    \pro(T^{\tra}) &=& \pro(T)^{\tra} \\
    \dem(T^{\tra}) &=& \dem(T)^{\tra} \\
    \pro(T^{\rev}) &=& \dem(T)^{\rev} \\
    \dem(T^{\rev}) &=& \pro(T)^{\rev}.
\end{eqnarray*}

Furthermore,
transposition and reversal act on the descent set of a tableau as follows.
For $i\in [n-1]$, $i\in \Des(T^{\tra})$ if and only if $i\not\in \Des(T)$,
and $i\in \Des(T^{\rev})$ if and only if $n-i \in \Des(T)$.

\begin{ex}
    Below is a tableau $T$ and its transpose and reverse.
    \begin{equation*}
        T = \young(::268,::3,14,59,7),
        \quad\quad
        T^{\tra} = \young(::157,::49,23,6,8),
        \quad\quad
        T^{\rev} = \young(::::3,:::15,:::69,::7,248).
    \end{equation*}
    Note that $\Des(T) = \{2,3,4,6,8\}$, $\Des(T^{\tra}) = \{1,5,7\}$,
    and $\Des(T^{\rev}) = \{1,3,5,6,7\}$,
    so the aforementioned descent relations hold.
\end{ex}

\subsection{Notation and Nomenclature Conventions}

For a number $x$ and standard Young tableau $T$,
we use $\pos_T(x)$ to denote the position of $x$ in $T$.

Where applicable, we will use lowercase letters to
denote the entries of standard Young tableaux
and uppercase letters to denote cells in standard Young tableaux.

We say that a cell $W$ is the \vocab{northern neighbor} of $Z$
if $W$ and $Z$ are adjacent and $W$ is north of $Z$.
We analogously define $W$ being the \vocab{eastern}, \vocab{southern},
and \vocab{western neighbors} of $Z$.

We say $W$ is \vocab{strictly north} of $Z$
if the row of $W$ is strictly north of the row of $Z$.
We say $W$ is \vocab{nonstrictly north} of $Z$ if $W$ is strictly north of $Z$
or in the same row as $Z$.
Note that a cell is nonstrictly north of itself,
and two unequal cells in the same row are nonstrictly north of each other.
We analogously define $W$ being (\vocab{strictly} or \vocab{nonstrictly})
\vocab{east}, \vocab{south}, or \vocab{west} of $Z$.

We say $W$ is \vocab{nonstrictly} (resp. \vocab{strictly}) \vocab{southwest}
of $Z$ if it is nonstrictly (resp. strictly) south and west of $Z$.
We analogously define $W$ being
\vocab{northwest}, \vocab{northeast}, and \vocab{southeast} of $Z$.

Finally, we say $W$ is \vocab{nonstrictly} (resp. \vocab{strictly})
\vocab{due north} of $Z$ if $W$ and $Z$ are in the same column
and $W$ is nonstrictly (resp. strictly) north of $Z$.
We analogously define $W$ being
\vocab{due east}, \vocab{south}, and \vocab{west} of $Z$.

Throughout this paper, if an assertion about relative positioning is made
without mention of strictness, the nonstrict form of that assertion is implied.

\section{Rotation}\label{sec:rotation}

In this section we define the rotation operation,
which will feature prominently in our main construction.
There are two variants of rotation:
southeast rotation, denoted $\RotSE$, and northwest rotation, denoted $\RotNW$.

These two operations are dual with respect to reversal,
in a sense that will be made precise later.
Thus we first define $\RotSE$, and state and prove its properties.

\subsection{Southeast Rotation}

\begin{figure}[H]
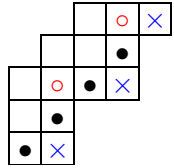
\label{fig:boundary}
    \centering
    \begin{equation*}
        \young(::\blank\rcirc\btimes,:\blank\blank\bull,\blank\rcirc\bull\btimes,\blank\bull,\bull\btimes)
    \end{equation*}
    \caption{
        The southeast boundary of a skew diagram.
        Exterior corners are marked with blue crosses $\btimes$,
        and interior corners are marked with red circles $\rcirc$.
        Other points on the southeast boundary are marked with black dots $\bull$.
    }
\end{figure}

The \vocab{southeast boundary} of the skew shape $\skewshape$ is
the set of cells in $\skewshape$ that do not have a diagonally-adjacent
southeastern neighbor.
Furthermore, a \vocab{southeast exterior corner} of $\skewshape$ is a cell
in the southeast boundary that has neither a southern nor eastern neighbor,
and a \vocab{southeast interior corner} is a cell in the southeast boundary
that has both southern and eastern neighbors.
Figure~\ref{fig:boundary} illustrates these definitions.

Within each connected component of $\skewshape$, the southeast boundary
is a connected ribbon, which has a natural southwest-to-northeast linear order.

A sequence of numbers $x_1,\ldots,x_k$ is
\vocab{min-unimodal} if there is $i\in [k]$ such that
\begin{equation*}
    x_1>x_2>\cdots>x_i<x_{i+1}<\cdots<x_k.
\end{equation*}

The notion of \vocab{southeast min-unimodality} will be central to our
construction.
\begin{defn}
    For $T\in \SYT(\skewshape)$ and a set $S\subseteq [n]$,
    $S$ is \vocab{southeast min-unimodal} in $T$
    if it satisfies the following two conditions.
    \begin{itemize}
        \item In $T$, all elements of $S$ lie in the same connected component
        of $\skewshape$ and are on the southeast boundary.
        \item The elements of $S$, in the southeast boundary's natural
        southwest-to-northeast order, form a min-unimodal sequence.
    \end{itemize}
\end{defn}

\begin{defn}
    Let $T\in \SYT(\skewshape)$.
    The \vocab{southeast rotation candidate set} of $T$, denoted $\RcSE(T)$,
    is the set $\{n,n-1,\ldots,n-k+1\}$, where $k$ is maximal such that
    this set is southeast min-unimodal in $T$.
    Let
    \begin{equation*}
        \RpSE(T) = \{\pos_T(x) | x\in \RcSE(T)\}
    \end{equation*}
    be the cells occupied by $\RcSE(T)$ in $T$.
\end{defn}

\begin{defn}
    The \vocab{southeast rotation endpoint} of $T$ is the southeast exterior
    corner due east or due south of $\pos_T(\min \RcSE(T))$.
    If $\pos_T(\min \RcSE(T))$ is a southeast interior corner,
    there are two such exterior corners;
    the southeast rotation endpoint is the one on the southeast boundary
    between $\pos_T(\min \RcSE(T))$ and $\pos_T(n)$.
\end{defn}

Note that $\pos_T(n)$ is always on the southeast boundary.
Also note that the southeast rotation endpoint is in $\RpSE(T)$
because its entry is at least $n-k+1$.

\begin{defn}
    Let $T\in \SYT(\skewshape)$.
    Let $x_1,x_2,\ldots,x_j$ denote the entries of $\RpSE(T)$
    between $\pos_T(n)$ and the southeast rotation endpoint of $T$, inclusive,
    in that order along the southeast boundary.
    Thus $x_1=n$.
    To construct $\RotSE(T)$ from $T$, move $n$ to $\pos_T(x_j)$,
    and for $i=1,\ldots,j-1$, move $x_{i+1}$ to $\pos_T(x_i)$.
\end{defn}

\begin{ex}\label{ex:rotse-move}
    Let
    \begin{equation*}
        T = \young(:279\btwelve\bthirteen,14\bten\beleven,36,58,\bfourteen,\bfifteen).
    \end{equation*}
    The set $\RcSE(T)$ is shown in {\bf bold}.
    The southeast rotation endpoint of $T$ is $\pos_T(11)$.
    So, the action of $\RotSE$ on $T$ is
    \begin{equation*}
        \young(:279\btwelve\bthirteen,14\rbten\rbeleven,36,58,\rbfourteen,\rbfifteen)
        \quad
        \overset{\RotSE}{\mapsto}
        \quad
        \young(:279\btwelve\bthirteen,14\rbeleven\rbfifteen,36,58,\rbten,\rbfourteen)
    \end{equation*}
    The entries rotated by $\RotSE$ are shown in \textcolor{red}{\bf red}.
\end{ex}

\begin{ex}
    Let
    \begin{equation*}
        T = \young(:1379\ten,2\beleven\bthirteen\bfifteen,4\btwelve,5\bfourteen,6,8)
    \end{equation*}
    The set $\RcSE(T)$ is shown in {\bf bold}.
    The southeast rotation endpoint of $T$ is $\pos_T(15)$.
    Since $\RotSE$ moves 15 to its own cell, $\RotSE(T) = T$.
\end{ex}

\begin{ex}
Any standard Young tableau $T$ of shape
\begin{equation*}
    \skewshape = \young(:::\blank,:\blank\blank\blank,\blank\blank\blank\blank)
\end{equation*}
has $n=8$ in the southeast corner.
Since the only southeast exterior corner of $\skewshape$
is the southeast corner,
the southeast rotation endpoint is that corner for all $T\in \SYT(\skewshape)$.
Therefore, $\RotSE$ is the identity map on $\SYT(\skewshape)$.
\end{ex}

\subsection{Properties of Southeast Rotation}

In this section, we prove some properties of southeast rotation.
For the sake of clarity, we begin by stating all the properties,
and then give their proofs.

Note that in a skew shape with more than one connected component,
the connected components have a strict southwest-to-northeast order,
in the sense that for any two distinct connected components,
all the cells in one are strictly northeast of all the cells of the other.

A \vocab{hook} is a connected shape consisting of a single row and a single
column joined at the northwest corner.
\begin{prop}\label{prop:rcse-locations}
    $\RpSE(T)$ is a skew shape
    contained in the southeast boundary of $\skewshape$.
    The southeast rotation endpoint $X$ and the cell $Y=\pos_T(\min \RcSE(T))$
    are in the same connected component $R$ of $\RpSE(T)$,
    which is one of the following.
    \begin{enumerate}[label=(\alph*)]
        \item A hook with northwestern corner $Y$ and southern or eastern endpoint $X$.
        $X$ is the southern endpoint if $\pos_T(n)$ is the southwesternmost cell in $\RpSE(T)$ and
        the eastern endpoint if $\pos_T(n)$ is the northeasternmost cell in $\RpSE(T)$.
        \item A single-column rectangle with northern and southern endpoints $Y$ and $X$, respectively.
        \item A single-row rectangle with western and eastern endpoints $Y$ and $X$, respectively.
    \end{enumerate}
    Moreover, all connected components of $\RpSE(T)$ strictly to the northeast
    of $R$ are single-row rectangles, and
    all connected components of $\RpSE(T)$ strictly to the southwest of $R$ are
    single-column rectangles.
\end{prop}
Note that a $1\times 1$ square is both a single-row rectangle and a
single-column rectangle.

This result motivates the following definition.
\begin{defn}
    Let $\RpSEc(T)$ denote the connected component of $\RpSE(T)$ containing
    $\pos_T(\min \RcSE(T))$.
    Let $\RpSEsw(T)$ and $\RpSEne(T)$ denote the subsets of $\RpSE(T)$ strictly
    southwest and strictly northeast of $\RpSEc(T)$, respectively.
\end{defn}
By Proposition~\ref{prop:rcse-locations},
$\RpSEc(T)$ is a hook, one-column rectangle, or one-row rectangle;
$\RpSEsw(T)$ is a disconnected union of one-column rectangles;
and, $\RpSEne(T)$ is a disconnected union of one-row rectangles.

\begin{ex}
    In the following four tableaux $T$, the entries in $\RpSE(T)$ are shown in {\bf bold},
    and the entries in $\RpSEc(T)$ are shown in \textcolor{red}{\bf red}.
    The first three tableaux, from left to right, correspond to
    cases (a), (b), and (c) of Proposition~\ref{prop:rcse-locations}.
    In the fourth, $\RpSEc(T)$ is a single cell, which corresponds to cases (b) and (c).
    \begin{eqnarray*}
        \young(:1247\bfourteen,3\rbnine\rbeleven\rbtwelve,5\rbten,6\rbthirteen,8,\bfifteen)
        \quad\quad
        \young(:236\btwelve\bthirteen,1578,4\rbten,9\rbeleven,\bfourteen,\bfifteen)
        \quad\quad
        \young(:279\btwelve\bthirteen,14\rbten\rbeleven,36,58,\bfourteen,\bfifteen)
        \quad\quad
        \young(:256\btwelve\bthirteen,147\rbnine,3\bten,8\beleven,\bfourteen,\bfifteen).
    \end{eqnarray*}
\end{ex}

\begin{prop}\label{prop:rotse-preserve-descents}
    For $T\in \SYT(\skewshape)$ and $i\in \{1,\ldots,n-2\}$,
    $i\in \Des(T)$ if and only if $i\in \Des(\RotSE(T))$.
\end{prop}

\begin{ex}
    For the $T$ in Example~\ref{ex:rotse-move}, we have
    $\Des(T) = \{2,4,7,9,13\}$ while $\Des(\RotSE(T)) = \{2,4,7,9\}$.
    The restrictions of the two descent sets to $\{1,\ldots,n-2\}$ agree.
\end{ex}

\begin{prop}\label{prop:rotse-automorphism}
    $\RotSE$ is a bijection from $\SYT(\skewshape)$ to itself.
\end{prop}

We will prove these properties in order, starting with
Proposition~\ref{prop:rcse-locations}.

\begin{proof}[Proof of Proposition~\ref{prop:rcse-locations}]
    $\RpSE(T)$ is, by definition, contained
    in the southeast boundary of $\skewshape$.
    Since $\RcSE(T)$ consists of the $|\RcSE(T)|$ largest elements
    of $\{1,\ldots,n\}$,
    $\RpSE(T)$ has the property that all cells of $\skewshape$ (nonstrictly)
    southeast of a cell in $\RpSE(T)$ are in $\RpSE(T)$.
    This implies that $\RpSE(T)$ is a skew shape.
    Since $X$ is due south or due east of $Y$,
    $X$ and $Y$ are in the same connected component of this skew shape.

    We claim that the shape formed by the cells
    \begin{equation*}
        R_{\NE} = \{Z\in \RpSE(T) | \text{$Z$ is nonstrictly northeast of $Y$}\}
    \end{equation*}
    avoids a vertical $2\times 1$ rectangle.

    This is because the entries in this subset of $\RpSE(T)$ is increasing from
    southwest to northeast.
    So, $R_{\NE}$ is a disconnected union of single-row rectangles.
    It follows that the connected components strictly northeast of
    the connected component $R$ are single-row rectangles;
    note that $R$ itself is not necessarily a single-row rectangle
    because it can have cells that are not northeast of $Y$.

    By an analogous argument, the subset
    \begin{equation*}
        R_{\SW} = \{Z\in \RpSE(T) | \text{$Z$ is nonstrictly southwest of $Y$}\}
    \end{equation*}
    is a disconnected union of single-column rectangles,
    so the connected components strictly southwest of $R$
    are single-column rectangles.

    If $Y$ is part of a single-row rectangle of length more than 1
    and a single-column rectangle of length more than 1, we get configuration (a).
    If $Y$ is part of a single-column rectangle, we get configuration (b).
    If $Y$ is part of a single-row rectangle, we get configuration (c).
\end{proof}

\begin{proof}[Proof of Proposition~\ref{prop:rotse-preserve-descents}]
    Let us assume $\pos_T(n)$ is the southwesternmost cell in $\RpSE(T)$.
    Let $\RcSE(T) = \{n,\ldots,n-k+1\}$.
    Let $X$ be the southeast rotation endpoint of $T$, and let $Y=\pos_T(n-k+1)$.
    Recall that $X\in \RpSEc(T)$.
    The action of $\RotSE$ takes one of two forms.
    \begin{enumerate}[label=(\alph*)]
        \item If $X$ is strictly due south of $Y$,
        we are in configurations (a) or (c) of Proposition~\ref{prop:rcse-locations}.
        Then, $\RotSE$ moves $n$ to $X$,
        each remaining entry of $\RpSEsw(T)$ to the next cell to the south
        in $\RpSEsw(T)$, and the entry in $X$ to the northernmost cell
        of $\RpSEsw(T)$.
        \item If $X$ is nonstrictly due east of $Y$,
        we are in configuration (b) of Proposition~\ref{prop:rcse-locations},
        and $\RpSEc(T)$ is a single-row rectangle.
        Then, $\RotSE$ moves $n$ to $X$,
        the remaining entries of $\RpSEsw(T)$ to the next cell to the south
        in $\RpSEsw(T)$, $n-k+1$ to the northernmost cell of $\RpSEsw(T)$,
        and the remaining entries of $\RpSEc(T)$ one cell to the west.
    \end{enumerate}
    We perform casework on $i$.

    Case 1: $i<n-k$.

    The action of $\RotSE$ affects neither $i$ nor $i+1$, so the lemma follows.

    Case 2: $i=n-k$.

    The action of $\RotSE$ does not move $n-k$.
    If it also does not move $n-k+1$, the lemma follows.
    Otherwise the action of $\RotSE$ must take the form (b).
    So, $\RpSEc(T)$ is a single-row rectangle,
    with western endpoint $Y$ and eastern endpoint $X$ (possibly with $X=Y$).

    Let $Z$ be the northernmost cell of $\RpSEsw(T)$.
    The action of $\RotSE$ moves $n-k+1$ from $Y$ to $Z$.
    Suppose for contradiction that $n-k$ is in one of $\Des(T)$ and $\Des(\RotSE(T))$ but not the other.
    Then $\pos_T(n-k)$ must be strictly north of $Z$ and nonstrictly south of $Y$.

    The cells nonstrictly due east of $Y$ are in $\RpSE(T)$.
    If $\pos_T(n-k)$ were strictly due west of $Y$, it must be the western neighbor of $Y$.
    Then, $\{n, \ldots, n-k\}$ is southeast min-unimodal in $T$, contradicting the maximality of $k$.
    Thus $\pos_T(n-k)$ is strictly south of $Y$.
    Since there are no cells in $\RpSE(T)$ strictly north of $Z$ and strictly south of $Y$,
    and $n-k$ is the largest entry of $T$ outside of $\RpSE(T)$,
    $n-k$ lies on the southeast boundary.
    But then, $\{n, \ldots, n-k\}$ is southeast min-unimodal in $T$, contradicting the maximality of $k$.

    Case 3: $i=n-k+1$.

    Note that for this case to occur, we must have $n-k+1 \le n-2$, so $n-k+2 \neq n$.

    If the action of $\RotSE$ takes the form (a),
    the set of cells whose entries are moved by $\RotSE$ is strictly south of $Y$.
    So, $\pos_T(n-k+1) = \pos_{\RotSE(T)}(n-k+1) = Y$,
    and $\pos_T(n-k+2)$ is strictly south of $Y$ if and only if $\pos_{\RotSE(T)}(n-k+2)$
    is strictly south of $Y$.
    Hence, $n-k+1$ is a descent of $T$ if and only if it is a descent of $\RotSE(T)$.

    If the action of $\RotSE$ takes the form (b),
    it moves $n-k+1$ from $Y$ to the northernmost cell of $\RpSEsw(T)$
    and $n$ from the southernmost cell of $\RpSEsw(T)$ to $X$.
    Aside from these entries, all entries in $\RpSEsw(T)$ remain in $\RpSEsw(T)$,
    and all entries in $\RpSEc(T) \cup \RpSEne(T)$ remain in $\RpSEc(T) \cup \RpSEne(T)$.

    Recall that $n-k+2\neq n$.
    Thus, either $\pos_T(n-k+2) \in \RpSEsw(T)$ and $\pos_{\RotSE(T)}(n-k+2) \in \RpSEsw(T)$,
    in which case $n-k+1$ is a descent of both $T$ and $\RotSE(T)$,
    or $\pos_T(n-k+2)\not\in \RpSEsw(T)$ and $\pos_{\RotSE(T)}(n-k+2) \not\in \RpSEsw(T)$,
    in which case $n-k+1$ is a descent of neither $T$ nor $\RotSE(T)$.

    Case 4: $n-k+1<i \le n-2$.

    By examining the actions (a) and (b), we see that for any $n-k+1<j_1,j_2\le n-1$,
    $\pos_T(j_1)$ is strictly south of $\pos_T(j_2)$ if and only if
    $\pos_{\RotSE(T)}(j_1)$ is strictly south of $\pos_{\RotSE(T)}(j_2)$.
    This implies that $i\in \Des(T)$ if and only if $i\in \Des(\RotSE(T))$.

    This proves the lemma when $\pos_T(n)$ is southwesternmost in $\RpSE(T)$.
    If $\pos_T(n)$ is northeasternmost in $\RpSE(T)$,
    then $\pos_{T^{\tra}}(n)$ is southwesternmost in $\RpSE(T^{\tra})$.
    The above argument shows that for $i\in [n-2]$,
    $i\in \Des(T^{\tra})$ if and only if $i\in \Des(\RotSE(T^{\tra}))$.
    But $i\in \Des(T)$ if and only if $i\not\in \Des(T^{\tra})$
    and $i\in \Des(\RotSE(T))$ if and only if $i\not\in \Des(\RotSE(T^{\tra}))$,
    where we use that $\RotSE$ commutes with transposition.
    Therefore, $i\in \Des(T)$ if and only if $i\in \Des(\RotSE(T))$.
\end{proof}

The rest of this subsection will be devoted to proving
Proposition~\ref{prop:rotse-automorphism}.

\begin{lem}\label{lem:rotse-image-syt}
    The image of $\RotSE$ is in $\SYT(\skewshape)$.
    That is, for any $T\in \SYT(\skewshape)$, the tableau $\RotSE(T)$ is standard.
\end{lem}
\begin{proof}
    We will show that the entries of $\RpSE(T)$ in $\RotSE(T)$ are standard.
    Since $\RotSE$ does not alter entries outside of $\RpSE(T)$,
    and all entries in $\RpSE(T)$ are larger than all entries not in
    $\RpSE(T)$, this proves that $\RotSE(T)$ is standard.

    By southeast min-unimodality of $\RotSE(T)$,
    $\pos_T(n)$ is the southwesternmost or northeasternmost cell in $\RpSE(T)$.
    Let us first assume it is southwesternmost.

    In $T$, let $\RcSE(T)$ be in the southwest-to-northeast order
    $x_1,\ldots,x_k$, with
    \begin{equation*}
        n=x_1 > \cdots > x_i = n-k+1 < x_{i+1} < \cdots < x_k,
    \end{equation*}
    and let $Y=\pos_T (n-k+1)$.
    Let $X$ be the southeast rotation endpoint of $T$,
    and let $x_j$ be the entry in $X$.

    $\RotSE$ does not move the entries in $\RpSEne(T)$,
    so the entries of $\RpSEne(T)$ are standard in $\RotSE(T)$.
    It remains to show that the entries of $\RpSEc(T)$ and $\RpSEsw(T)$
    are standard in $\RotSE(T)$.

    We consider two cases.

    Case 1: $X$ is strictly due south of $Y$.

    In this case $j<i$, and the set of rotated entries does not contain $x_i$.
    This case corresponds to configuration (a) of Proposition~\ref{prop:rcse-locations},
    or configuration (b) where $X\neq Y$.

    $\RpSEsw(T)$ is a union of one-column rectangles, arranged from southwest to northeast.
    In $T$, the entries of $\RpSEsw(T)$ are, from south to north, $x_1,\ldots,x_{j-1}$.
    Thus, in $\RotSE(T)$, the entries of $\RpSEsw(T)$ are, from south to north, $x_2,\ldots,x_j$.
    This is standard.

    $\RotSE$ moves $n$ to $X$ but does not otherwise change the entries in $\RpSEc(T)$.
    Since $X$ is the southern endpoint of either a hook or a single-column rectangle,
    the entries of $\RpSEc(T)$ are standard in $\RotSE(T)$.

    Case 2: $X$ is nonstrictly due east of $Y$.

    In this case $j\ge i$, and the set of rotated elements contains $x_i$.
    This case corresponds to configuration (c) of Proposition~\ref{prop:rcse-locations},
    where possibly $X=Y$.

    In $T$, the entries of $\RpSEsw(T)$ are, from south to north, $x_1,\ldots,x_{i-1}$,
    and the entries of $\RpSEc(T)$ are, from east to west, $x_j,x_{j-1},\ldots,x_i$.
    So, in $\RotSE(T)$, the entries of $\RpSEsw(T)$ are, from south to north, $x_2,\ldots,x_i$,
    and the entries of $\RpSEc(T)$ are, from east to west, $x_1,x_j,x_{j-1},\ldots,x_{i+1}$.
    Both sequences are decreasing, so the entries of $\RpSEc(T)$ and $\RpSEsw(T)$
    are standard in $\RotSE(T)$.

    This proves the lemma when $\pos_T(n)$ is southwesternmost in $\RpSE(T)$.
    If $\pos_T(n)$ is northeasternmost in $\RpSE(T)$,
    then $\pos_{T^{\tra}}(n)$ is southwesternmost in $\RpSE(T^{\tra})$.
    The above argument shows that $\RotSE(T^{\tra})$ is standard.
    So, $\RotSE(T) = \RotSE(T^{\tra})^{\tra}$ is standard.
\end{proof}

For $S\subset [n]$, let
\begin{equation*}
    \pos_T(S) = \{\pos_T(x) | x\in S\}
\end{equation*}
denote the cells in $T$ occupied by $S$.

To construct the inverse of $\RotSE$, we will need a notion of \vocab{balance}.

\begin{defn}
    If $S$ is a southeast min-unimodal set in $T$,
    and $W$ is an arbitrary cell in the same connected component as $\pos_T(S)$
    and on the southeast boundary,
    $S$ is \vocab{balanced} with respect to $W$ in $T$ if the entries of
    \begin{equation*}
        \pos_T(S)_W^{\NE} =
        \{W'\in \pos_T(S) | \text{$W'$ is nonstrictly northeast of $W$}\}
    \end{equation*}
    and
    \begin{equation*}
        \pos_T(S)_Z^{\SW} =
        \{W'\in \pos_T(S) | \text{$W'$ is nonstrictly southwest of $W$}\}
    \end{equation*}
    are both increasing in the direction away from $W$.
\end{defn}
Equivalently, since $S$ is southeast min-unimodal,
$S$ is balanced with respect to $W$ if $W$ is consecutive with $\pos_T(\min S)$
in the southwest-to-northeast order of $\pos_T(S) \cup W$
(or equal to $\pos_T(\min S)$).

\begin{defn}
    The \vocab{northern balance point} of a southeast exterior corner $Z$
    is the farthest point on the southeast boundary due north of $Z$.
    The \vocab{western balance point} of $Z$ is defined analogously.
\end{defn}
Equivalently, if $Z$ is the northeasternmost southeast exterior corner
in its connected component of $\skewshape$,
its northern balance point is the northeastern corner of its connected component,
and otherwise it is the southeast interior corner due north of $Z$.
The western balance point of $Z$ is characterized analogously.

The following lemma allows us to recover $\RcSE(T)$
from the tableau $\RotSE(T)$.
\begin{lem}\label{lem:recover-rcse}
    The set $\RcSE(T)$ is $\{n,\ldots,n-k+1\}$,
    where $k$ is maximal such that $S=\{n-1,\ldots,n-k+1\}$
    obeys the following two properties.
    \begin{enumerate}[label = (\roman*)]
        \item In $R=\RotSE(T)$, $S$ is southeast min-unimodal and contained
        in the same connected component as $n$.
        \item The set $\pos_R(S)$ does not contain
        both the northern and western neighbors of $\pos_R(n)$.
        Moreover:
        \begin{enumerate}
            \item If $\pos_R(S)$ contains the northern (resp. western) neighbor of $\pos_R(n)$,
            $S$ is balanced with respect to the northern (resp. western) balance point of $\pos_R(n)$.
            \item If $\pos_R(S)$ contains neither neighbor of $\pos_R(n)$,
            $S$ is balanced with respect to both the northern and western balance points of $\pos_R(n)$.
        \end{enumerate}
    \end{enumerate}
\end{lem}
Note that properties (i) and (ii) become more restrictive as $k$ increases.
So, for each $T$, $S=\{n-1,\ldots,n-k+1\}$ satisfies these properties
for all $k$ up to a threshold value, and for no $k$ greater than it;
the maximal $k$ equals the threshold value.
\begin{proof}
    Let $\RcSE(T) = \{n,\ldots,n-k'+1\}$,
    and let its elements appear in the southwest-to-northeast order
    $x_1,\ldots,x_{k'}$, such that
    \begin{equation*}
        x_1 > \cdots > x_i = n-k'+1 < x_{i+1} < \cdots < x_{k'}.
    \end{equation*}
    Let $S' = \{n-1,\ldots,n-k'+1\}$.  We will show that $S' = S$.

    Let $X$ be the southeast rotation endpoint of $T$, so $X=\pos_{R}(n)$.
    Furthermore, let $Y = \pos_T(n-k'+1)$.
    Since $\RotSE$ preserves the southwest-to-northeast order of $S'$,
    this set is southeast min-unimodal in $R$.
    Since $\RcSE(T)$ is contained in the same connected component of $T$,
    all of $S'$ is in the same connected component of $R$ as $n$.
    Thus, $S'$ satisfies condition (i).

    By Proposition~\ref{prop:rcse-locations},
    $\RpSE(T)$ cannot contain both the western and northern neighbors of $X$.
    By examining the three configurations in
    Proposition~\ref{prop:rcse-locations},
    we see that if $X$ is strictly due south (resp. strictly due east) of $Y$,
    $S'$ is balanced with respect to the northern (resp. western) balance point of $X$,
    and if $X=Y$, $S'$ is balanced with respect to both balance points of $X$.
    Thus, $S'$ satisfies condition (ii).

    It remains to show that $S''=\{n-1,\ldots,n-k'\}$
    does not satisfy both conditions (i) and (ii).
    If $n-k'=0$ there is nothing to prove,
    so we can assume $n-k'\ge 1$.
    Suppose for contradiction that $S''$ satisfies conditions (i) and (ii).
    We consider three cases.

    Case 1: $X$ is strictly due south of $Y$.

    As $\pos_R(S')$ contains the northern neighbor of $X$, so does $\pos_{R}(S'')$.
    So, $S''$ is balanced in $R$ with respect to the northern balance point $X_N$ of $X$.
    If $Y=X_N$, $Y$ is either a southeast interior corner or the northeastern corner
    of a connected component of $\lambda/\mu$.
    In either case, $\pos_T(n)$ must be southwest of $Y$,
    and thus (by definition of southeast rotation endpoint) southwest of $X$.
    But then, $\RotSE$ does not move the entry in $Y$, so $\pos_R(n-k'+1) = Y = X_N$.
    This makes it impossible for both $\pos_R (S'')_Y^{\NE}$ and $\pos_R (S'')_Y^{\SW}$
    to be increasing away from $X_N$, because one of these sets of cells contains the entry $n-k'$.
    Thus $Y\neq X_N$.

    The cell in $\RpSE(T)$ southwest of $X_N$ closest to $X_N$ is $Y=\pos_T(x_i)$.
    If $i<k'$, the cell in $\RpSE(T)$ northeast of $X_N$ closest to $X_N$ is $\pos_T(x_{i+1})$;
    since $S''$ is balanced in $R$ with respect to $X_N$,
    $\pos_{R}(n-k')$ is on the southeast boundary between $\pos_T(x_i)$ and $\pos_T(x_{i+1})$.
    If $i=k'$, there are no cells in $\RpSE(T)$ northeast of $X_N$;
    since $S''$ is balanced in $R$ with respect to $X_N$,
    $\pos_{R}(n-k')$ is on the southeast boundary northeast of $\pos_T(x_i)$.
    In either case, $\RotSE$ does not move $n-k'$, so $\pos_{R}(n-k') = \pos_T(n-k')$.
    So, $\{n,\ldots,n-k'\}$ is southeast min-unimodal in $T$.
    This contradicts the maximality of $k'$ in the definition of $\RcSE(T)$.

    Case 2: $X$ is strictly due east of $Y$.

    The argument here is symmetric.

    As $\pos_R(S')$ contains the northern neighbor of $X$, so does $\pos_{R}(S'')$.
    So, $S''$ is balanced in $R$ with respect to the western balance point $X_W$ of $X$.
    If $Y=X_W$, $Y$ is either a southeast interior corner or the southwestern corner
    of a connected component of $\lambda/\mu$.
    In either case, $\pos_T(n)$ must be northeast of $Y$,
    and thus (by definition of southeast rotation endpoint) northeast of $X$.
    But then, $\RotSE$ does not move the entry in $Y$, so $\pos_R(n-k'+1) = Y = X_W$.
    This makes it impossible for both $\pos_R (S'')_Y^{\NE}$ and $\pos_R (S'')_Y^{\SW}$
    to be increasing away from $X_W$, because one of these sets of cells contains the entry $n-k'$.
    Thus $Y\neq X_W$.

    The cell in $\RpSE(T)$ northeast of $X_W$ closest to $X_W$ is $Y=\pos_T(x_i)$.
    If $i>1$, the cell in $\RpSE(T)$ southwest of $X_W$ closest to $X_W$ is $\pos_T(x_{i-1})$;
    since $S''$ is balanced in $R$ with respect to $X_W$,
    $\pos_{R}(n-k')$ is on the southeast boundary between $\pos_T(x_i)$ and $\pos_T(x_{i-1})$.
    If $i=1$, there are no cells in $\RpSE(T)$ southwest of $X_W$;
    since $S''$ is balanced in $R$ with respect to $X_W$,
    $\pos_{R}(n-k')$ is on the southeast boundary southwest of $\pos_T(x_i)$.
    In either case, $\RotSE$ does not move $n-k'$, so $\pos_{R}(n-k') = \pos_T(n-k')$.
    So, $\{n,\ldots,n-k'\}$ is southeast min-unimodal in $T$.
    This contradicts the maximality of $k'$ in the definition of $\RcSE(T)$.

    Case 3: $X=Y$.

    In this case, $\RpSEc(T)$ is a $1\times 1$ square consisting of $X$.
    Since all southern and eastern neighbors of cells in $\RpSE(T)$ are in $\RpSE(T)$,
    no cells northwest of $X$ are in $\RpSE(T)$ except $X$ itself.
    In particular, the balance points $X_N$ and $X_W$ are not in $\RpSE(T)$.
    Note that $S''$ is balanced in $R$ with respect to at least one balance point $Z\in \{X_N, X_W\}$.

    Let us first address the case $k'=1$.
    If $k'=1$, $\{n-1, n\}$ is not southeast min-unimodal in $T$,
    so $n-1$ is in a different connected component than $n$ in $T$.
    Thus the action of $\RotSE$ on $T$ is the identity map, and $R=T$.
    Then, $S'' = \{n-1\}$ does not satisfy condition (i) because it is not
    contained in the same connected component as $n$ in $R$.

    Otherwise, let us assume $k'>1$.
    If $1<i<k'$, the cell in $\RpSE(T)\setminus \{X\}$ southwest of $Z$ closest to $Z$ is $\pos_T(x_{i-1})$,
    and the cell in $\RpSE(T)\setminus \{X\}$ northeast of $Z$ closest to $Z$ is $\pos_T(x_{i+1})$;
    since $S''$ is balanced in $R$ with respect to $Z$,
    $\pos_{R}(n-k')$ is on the southeast boundary between $\pos_T(x_{i-1})$ and $\pos_T(x_{i+1})$.
    If $i=k$ (resp. $i=1$), the cell in $\RpSE(T)\setminus \{X\}$ southwest (resp. northeast) of $Z$
    closest to $Z$ is $\pos_T(x_{i-1})$ (resp. $\pos_T(x_{i+1})$),
    and there are no cells in $\RpSE(T)\setminus \{X\}$ northeast (resp. southwest) of $Z$;
    since $S''$ is balanced in $R$ with respect to $Z$,
    $\pos_{R}(n-k')$ is on the southeast boundary northeast (resp. southwest) of $Z$.
    In any case, $\RotSE$ does not move $n-k'$, so $\pos_{R}(n-k') = \pos_T(n-k')$.
    So, $\{n,\ldots,n-k'\}$ is southeast min-unimodal in $T$.
    This contradicts the maximality of $k'$ in the definition of $\RcSE(T)$.
\end{proof}

We now construct a function
$\Theta : \RotSE(\SYT(\skewshape)) \functo \SYT(\skewshape)$,
which we will show is the inverse of $\RotSE$
(and thus that $\RotSE(\SYT(\skewshape)) = \SYT(\skewshape)$).

Given $R=\RotSE(T)$, we can, using Lemma~\ref{lem:recover-rcse}, recover $\RcSE(T)$.
Since $\RotSE$ permutes the entries of $\RcSE(T)$ amongst themselves,
we also recover $\RpSE(T)$ as the set of positions of $\RcSE(T)$ in $R$.

$\Theta$ moves $n$ from $\pos_R(n)$ to either
the northeasternmost or southwesternmost position in $\RpSE(T)$
and displaces the entries of $\RpSE(T)$ in between one position closer to $\pos_R(n)$.
The direction it moves $n$ is as follows.
If $\pos_R(n)$ and $\pos_R(\min \RpSE(T))$ are in different connected components of $\RpSE(T)$,
$\Theta$ moves $n$ in the direction of $\pos_R(\min \RpSE(T))$;
otherwise, it moves $n$ in the opposite direction.

\begin{ex}
    Let
    \begin{equation*}
        R=\young(:279\btwelve\bthirteen,14\beleven\bfifteen,36,58,\bten,\bfourteen).
    \end{equation*}
    The set $\RcSE(T)$, as determined by Lemma~\ref{lem:recover-rcse},
    is in {\bf bold}.
    As $\pos_R(10)$ is southwest of $\pos_R(15)$
    and not in the same connected component,
    $\Theta$ moves $15$ to the southwesternmost cell of $\RpSE(T)$.
    Thus the action of $\Theta$ is
    \begin{equation*}
        \young(:279\btwelve\bthirteen,14\rbeleven\rbfifteen,36,58,\rbten,\rbfourteen)
        \quad
        \overset{\Theta}{\mapsto}
        \quad
        \young(:279\btwelve\bthirteen,14\rbten\rbeleven,36,58,\rbfourteen,\rbfifteen).
    \end{equation*}
    The entries moved by $\Theta$ are in \textcolor{red}{\bf red}.
\end{ex}
\begin{ex}
    Let
    \begin{equation*}
        R=\young(:279\btwelve\bthirteen,14\bten\bfifteen,36,58,\beleven,\bfourteen).
    \end{equation*}
    The set $\RcSE(T)$, as determined by Lemma~\ref{lem:recover-rcse},
    is in {\bf bold}.
    As $\pos_R(10)$ is southwest of $\pos_R(15)$
    but in the same connected component,
    $\Theta$ moves $15$ to the northeasternmost cell of $\RpSE(T)$.
    Thus, the action of $\Theta$ is
    \begin{equation*}
        \young(:279\rbtwelve\rbthirteen,14\bten\rbfifteen,36,58,\beleven,\bfourteen)
        \quad
        \overset{\Theta}{\mapsto}
        \quad
        \young(:279\rbthirteen\rbfifteen,14\bten\rbtwelve,36,58,\beleven,\bfourteen).
    \end{equation*}
    The entries moved by $\Theta$ are in \textcolor{red}{\bf red}.
\end{ex}
It isn't difficult to check that, in these examples, $\RotSE(\Theta(R))=R$.

\begin{lem}\label{lem:rotse-invertible}
    $\Theta$ is the left inverse of $\RotSE$.
    That is, $\Theta(\RotSE(T)) = T$ for all $T\in \SYT(\skewshape)$.
\end{lem}
\begin{proof}
    Let us assume $\pos_T(n)$ is the southwesternmost cell in $\RpSE(T)$.
    The case where $\pos_T(n)$ is northeasternmost is analogous.

    As before, let $X$ be the southeast rotation endpoint of $T$.
    Moreover, let $Y=\pos_T(\min \RcSE(T))$ and $R=\RotSE(T)$.
    We will first show that the action of $\Theta$ on $R$ sends $n$
    back to its original location.
    We consider two cases.

    Case 1: $X$ is strictly due south of $Y$.

    In this case, $\min \RcSE(T)$ is not among the rotated elements.
    Thus, $\pos_R(n)=X$ and $\pos_R(\min \RcSE(T))=Y$.
    Since $\pos_R(\min \RcSE(T))$ is northeast of $\pos_R(n)$ and in the same
    connected component,
    the action of $\Theta$ on $R$
    sends $n$ back to the southwesternmost cell in $\RpSE(T)$.

    Case 2: $X$ is nonstrictly due east of $Y$.

    In this case, $\min \RcSE(T)$ is among the rotated elements.
    Thus, $\pos_R(n)=X$ and $\pos_R(\min \RcSE(T))$
    is the northernmost cell in $\RpSEsw(T)$.
    Since $\pos_R(\min \RcSE(T))$ is southwest of $\pos_R(n)$
    and in a different connected component,
    the action of $\Theta$ on $R$
    sends $n$ back to the southwesternmost cell in $\RpSE(T)$.

    In either case, since $\RotSE$ and $\Theta$ preserve
    the southwest-to-northeast ordering of $\RcSE(T)\setminus \{n\}$,
    $\Theta$ also sends the remaining elements of $\RcSE(T)$
    to their original locations.
\end{proof}

\begin{proof}[Proof of Proposition~\ref{prop:rotse-automorphism}]
    By Lemma~\ref{lem:rotse-invertible}, the map $\RotSE$ is invertible on the left,
    and thus injective.
    By Lemma~\ref{lem:rotse-image-syt},
    the image $\RotSE(\SYT(\skewshape))$ is contained in $\SYT(\skewshape)$.
    Thus the domain and codomain of $\RotSE$ have the same size;
    this implies that $\RotSE$ is a bijective map on $\SYT(\skewshape)$.
\end{proof}
\subsection{Northwest Rotation}

Define the \vocab{northwest boundary}, \vocab{northwest exterior corners},
and \vocab{northwest interior corners} of $\skewshape$ analogously as above.
Within each connected component of $\skewshape$, the northwest boundary is
a connected ribbon with a natural southwest-to-northeast linear order.

A sequence of numbers $x_1,\ldots,x_k$ is
\vocab{max-unimodal} if there is $i\in [k]$ such that
\begin{equation*}
    x_1<x_2<\cdots<x_i>x_{i+1}>\cdots>x_k.
\end{equation*}

\begin{defn}
    For $T\in \SYT(\skewshape)$, a set $S\subseteq [n]$
    is \vocab{northwest max-unimodal} in $T$
    if it satisfies the following two conditions.
    \begin{itemize}
        \item In $T$, all elements of $S$ lie in the same connected component
        of $\skewshape$ and are on the northwest boundary.
        \item The elements of $S$, in the northwest boundary's natural
        southwest-to-northeast order, form a max-unimodal sequence.
    \end{itemize}
\end{defn}

The definitions below are analogous to those given above.

\begin{defn}
    Let $T\in \SYT(\skewshape)$.
    The \vocab{northwest rotation candidate set} of $T$, denoted $\RcNW(T)$,
    is the set $\{1,2,\ldots,\ell\}$, where $\ell$ is maximal such that
    this set is northwest max-unimodal in $T$.
    Let
    \begin{equation*}
        \RpNW(T) = \{\pos_T(x) | x\in \RcNW(T)\}
    \end{equation*}
    be the cells occupied by $\RcNW(T)$ in $T$.
\end{defn}

\begin{defn}
    The \vocab{northwest rotation endpoint} of $T$ is the northwest exterior
    corner due west or due north of $\pos_T(\max \RcNW(T))$.
    If $\pos_T(\max \RcNW(T))$ is a northwest interior corner,
    there are two such exterior corners;
    the northwest rotation endpoint is the one on the northwest boundary
    between $\pos_T(\max \RcNW(T))$ and $\pos_T(1)$.
\end{defn}

\begin{defn}
    Let $T\in \SYT(\skewshape)$.
    Let $x_1,x_2,\ldots,x_j$ denote the entries of $\RpNW(T)$
    between $\pos_T(1)$ and the northwest rotation endpoint of $T$, inclusive,
    in that order along the northwest boundary.
    Thus $x_1=1$.
    To construct $\RotNW(T)$ from $T$, move $1$ to $\pos_T(x_j)$,
    and for $i=1,\ldots,j-1$, move $x_{i+1}$ to $\pos_T(x_i)$.
\end{defn}

\begin{rem}
The two rotation operators commute with transposition, i.e.,
\begin{eqnarray*}
\RotSE(T^{\tra}) &=& \RotSE(T)^{\tra} \\
\RotNW(T^{\tra}) &=& \RotNW(T)^{\tra}
\end{eqnarray*}
for all $T\in \SYT(\skewshape)$.
Moreover, they are conjugate with respect to reversal, i.e.
\begin{eqnarray*}
\RotSE(T^{\rev}) &=& \RotNW(T)^{\rev} \\
\RotNW(T^{\rev}) &=& \RotSE(T)^{\rev}.
\end{eqnarray*}
for all $T\in \SYT(\skewshape)$.
\end{rem}

\subsection{Properties of Northwest Rotation}

The following analogues of Propositions~\ref{prop:rcse-locations}.
\ref{prop:rotse-preserve-descents}, and \ref{prop:rotse-automorphism}
hold for northwest rotation.
They follow from similar arguments.

A \vocab{reverse hook} is a connected shape consisting of a single row and a
single column joined at the southeast corner.
\begin{prop}\label{prop:rcnw-locations}
    $\RpNW(T)$ is a skew shape
    contained in the northwest boundary of $\skewshape$.
    The northwest rotation endpoint $X$ and the cell $Y=\pos_T(\max \RcNW(T))$
    are in the same connected component $R$ of $\RpSE(T)$,
    which is one of the following.
    \begin{enumerate}[label=(\alph*)]
        \item A reverse hook with southeastern corner $Y$
        and northern or western endpoint $X$.
        $X$ is the northern endpoint if $\pos_T(1)$ is the northeasternmost cell
        in $\RpNW(T)$ and
        the western endpoint if $\pos_T(1)$ is the southwesternmost cell
        in $\RpNW(T)$.
        \item A single-column rectangle with southern and northern
        endpoints $Y$ and $X$, respectively.
        \item A single-row rectangle with eastern and western
        endpoints $Y$ and $X$, respectively.
    \end{enumerate}
    Moreover, all connected components of $\RpNW(T)$ strictly to the northeast
    of $R$ are single-column rectangles, and
    all connected components of $\RpNW(T)$ strictly to the southwest of $R$ are
    single-row rectangles.
\end{prop}
\begin{prop}\label{prop:rotnw-preserve-descents}
For $T\in \SYT(\skewshape)$ and $i\in \{2,\ldots,n-1\}$,
$i\in \Des(T)$ if and only if $i\in \Des(\RotNW(T))$.
\end{prop}

\begin{prop}\label{prop:rotnw-automorphism}
$\RotNW$ is a bijection from $\SYT(\skewshape)$ to itself.
\end{prop}

\subsection{A Non-Interference Lemma}

Finally, we will need that these two rotation operators do not interfere with
each other, in the sense of Lemma~\ref{lem:rot-no-interference} below.
Critically, Lemma~\ref{lem:rot-no-interference} uses the hypothesis
that $\skewshape$ is not a connected ribbon.
On connected ribbons, this non-interference property can fail;
this part of our argument motivates the need for $\skewshape$ to be not a connected ribbon.

\begin{lem}\label{lem:rot-no-interference-aux}
    Suppose $T\in \SYT(\skewshape)$,
    and neither $\RcSE(T)$ nor $\RcNW(T)$ contains some $x\in [n]$.
    Then $\RcSE(T)$ and $\RcNW(T)$ are disjoint.
    Moreover, $\RcSE(T)=\RcSE(\RotNW(T))$ and $\RcNW(T)=\RcNW(\RotSE(T))$.
\end{lem}
\begin{proof}
    Since all entries of $\RcNW(T)$ are smaller than $x$,
    and all entries of $\RcSE(T)$ are larger than $x$, these sets are disjoint.
    Moreover, since $x$ is in neither set, $\min \RcSE(T)-1 \not\in \RcNW(T)$.
    Therefore, applying $\RotNW$ to $T$ does not move any elements of $\RcSE(T)$,
    nor does it move $\min \RcSE(T)-1$.

    In $T$, $\RcSE(T)$ is southeast min-unimodal, but $\RcSE(T)\cup \{\min \RcSE(T)-1\}$ is not;
    since none of the elements of $\RcSE(T)\cup \{\min \RcSE(T)-1\}$
    are moved by the action of $\RotNW$ on $T$, the same is true of $\RotNW(T)$.
    Thus, $\RcSE(T)=\RcSE(\RotNW(T))$.
    An analogous argument shows $\RcNW(T)=\RcNW(\RotSE(T))$.
\end{proof}

\begin{lem}\label{lem:rot-no-interference-aux2}
    Suppose $\skewshape$ contains a $2\times 2$ square.
    For any $T\in \SYT(\skewshape)$, each of $\RpSE(T)$ and $\RpNW(T)$
    contains at most two cells of this square.
\end{lem}
\begin{proof}
    Suppose for contradiction that $\RpSE(T)$ contains at least three cells of the square.
    Since all elements of $\RcSE(T)$ are larger than all elements of $[n]\setminus \RcSE(T)$,
    $\RpSE(T)$ contains the northeast, southeast, and southwest cells of the square.
    Thus, $\RpSE(T)$ contains the shape
    \begin{equation*}
        \young(:\blank,\blank\blank).
    \end{equation*}
    But Proposition~\ref{prop:rcse-locations} implies that $\RpSE(T)$ avoids this shape, contradiction.
    We get a similar contradiction if $\RpNW(T)$ contains at least three cells of the square.
\end{proof}

\begin{lem}\label{lem:rot-no-interference}
    Let $T\in \SYT(\skewshape)$, where $\skewshape$ is not a connected ribbon.
    Then the sets $\RcSE(T)$ and $\RcNW(T)$ are disjoint.
    Moreover, $\RcSE(T)=\RcSE(\RotNW(T))$ and $\RcNW(T)=\RcNW(\RotSE(T))$.
\end{lem}
\begin{proof}
    If $\skewshape$ is disconnected
    and $\pos_T(1)$ and $\pos_T(n)$ are in different connected components,
    then $\RpSE(T)$ and $\RpNW(T)$ are in different connected components.
    Hence, $\RcSE(T)$ and $\RcNW(T)$ are disjoint.
    Moreover,
    $\RotSE$ does not affect cells in the connected component of $\RpNW(T)$
    and vice versa;
    thus, we have $\RcNW(T)=\RcNW(\RotSE(T))$ and $\RcSE(T)=\RcSE(\RotNW(T))$.

    If $\skewshape$ is disconnected
    and $\pos_T(1)$ and $\pos_T(n)$ are in the same connected component,
    then there exists some $x\in [n]$ such that $\pos_T(x)$
    is not in this connected component.
    This $x$ is in neither $\RcSE(T)$ nor $\RcNW(T)$,
    so the lemma follows from Lemma~\ref{lem:rot-no-interference-aux}.

    Otherwise, $\skewshape$ is connected.
    By Lemma~\ref{lem:rot-no-interference-aux}, we can assume
    \begin{equation*}
        \RcSE(T)\cup \RcNW(T) = [n].
    \end{equation*}
    As $\skewshape$ is not a ribbon, it contains a $2\times 2$ square.
    Suppose that in $T$, the entries of this square are
    \begin{equation*}
        \young(uv,wz),
    \end{equation*}
    and label the corresponding positions $U, V, W, Z$.

    Each of $u,v,w,z$ belongs to $\RcSE(T)$ or $\RcNW(T)$.
    By Lemma~\ref{lem:rot-no-interference-aux2}, each of $\RcSE(T)$ and
    $\RcNW(T)$ contains at most two of $u,v,w,z$;
    as $\RcSE(T)\cup \RcNW(T) = [n]$, exactly two are in each set.
    Thus, either $u,w\in \RcNW(T)$ and $v,z \in \RcSE(T)$, or $u,v\in \RcNW(T)$
    and $w,z\in \RpSE(T)$.
    Suppose first that $u,w\in \RcNW(T)$ and $v,z \in \RcSE(T)$.

    Suppose for contradiction that $\RcSE(T)$ and $\RcNW(T)$ are not disjoint,
    so there is $a\in \RcSE(T)\cap \RcNW(T)$.
    Since $a\in \RcSE(T)$ while $w\not\in \RcSE(T)$, $a>w$.
    By northwest max-unimodality of $\RcNW(T)$,
    the entries of $\RpNW(T)$ are decreasing northeast of $W$.
    Thus, $\pos_T(a)$ is (nonstrictly) southwest of $W$.
    Since $a\in \RcNW(T)$ while $v\not\in \RcNW(T)$, $a<v$.
    By southeast min-unimodality of $\RcSE(T)$,
    the entries of $\RpSE(T)$ are increasing southwest of $V$.
    Thus, $\pos_T(a)$ is (nonstrictly) northeast of $V$.
    This is a contradiction. Thus, $\RcSE(T)$ and $\RcNW(T)$ are disjoint.

    Let $x = \max \RcNW(T)$.
    By northeast max-unimodality of $\RcNW(T)$, $\pos_T(x)$ is southwest of $W$.
    The rotation operation can only bring $x$ northeast as far as
    the next-northeastern position in $\RpNW(T)$, which is southwest of $U$.
    Thus, $\pos_{\RotNW(T)}(x)$ is southwest of $U$.

    Let $y = \min \RcSE(T)$; note that $y=x+1$.
    Because $\RcSE(T)$ is southeast min-unimodal, $\pos_T(y)$ is northeast of $V$.
    Since $\RcSE(T)$ and $\RcNW(T)$ are disjoint, $\RcSE(T) \subseteq \RcSE(\RotNW(T))$.
    If $\RcSE(\RotNW(T))$ and $\RcSE(T)$ are different, namely if $x\in \RcSE(\RotNW(T))$,
    by southeast min-unimodality $\pos_{\RotNW(T)}(x)$ must be a position on the southeast boundary
    between $\pos_T(y)$ and the next-northeastern or next-southwestern cell in $\RpSE(T)$.
    Thus, $\pos_{\RotNW(T)}(x)$ must be northeast of $V$, contradiction.
    Thus, $\RcSE(T)=\RcSE(\RotNW(T))$, and analogously, $\RcNW(T)=\RcNW(\RotSE(T))$.

    The case where $u,v\in \RcNW(T)$ and $w,z\in \RcSE(T)$ is analogous.
\end{proof}

\begin{ex}
    When $\skewshape$ is a connected ribbon, Lemma~\ref{lem:rot-no-interference} does not hold.
    When
    \begin{equation*}
        T = \young(:::2,:367,:4,15),
    \end{equation*}
    the sets $\RcNW(T) = \{1,2,3,4,5\}$ and $\RcSE(T) = \{3,4,5,6,7\}$ overlap.
\end{ex}

\section{Discussion of Theorem~\ref{thm:main}}\label{sec:discussion}

Recall our main result, reproduced below for clarity.
\begin{repthm}{thm:main}
Suppose $\skewshape$ is a skew shape that is not a connected ribbon.
Let
\begin{equation*}
    \phi = \RotNW^{-1} \circ \pro \circ \RotSE,
\end{equation*}
and let $n\in \cDes(T)$ if and only if $1\in \Des(\phi T)$.
Then $(\cDes, \phi)$ is a cyclic descent map.
\end{repthm}

We first present some examples of this result.

\begin{ex}\label{ex:phi}
    The action of $\phi$ on
    \begin{equation*}
        T = \young(:24,:35,16)
    \end{equation*}
    is as follows.
    \begin{equation*}
        \young(:24,:35,16)
        \quad
        \overset{\RotSE}{\mapsto}
        \quad
        \young(:24,:36,15)
        \quad
        \overset{\pro}{\mapsto}
        \quad
        \young(:13,:45,26)
        \quad
        \overset{\RotNW^{-1}}{\mapsto}
        \quad
        \young(:23,:45,16)
        .
    \end{equation*}
\end{ex}

\begin{ex}
Below are the orbits of $\phi$ on $\SYT((3,3,2)/(1,1))$,
along with the corresponding cyclic descents.

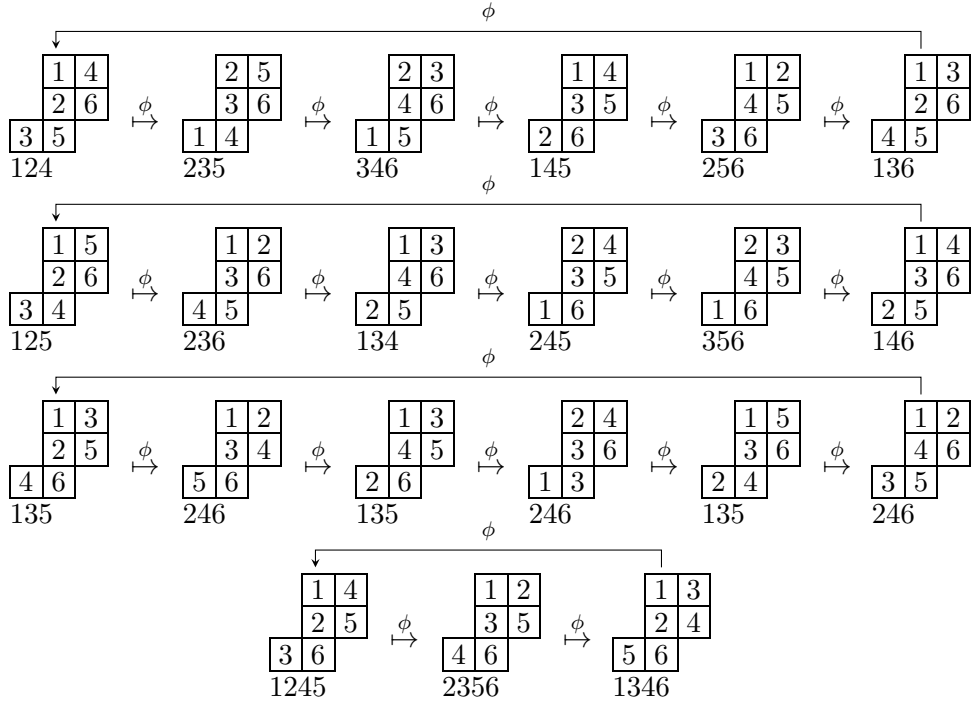
\begin{tikzpicture}[scale=2.3]
\node (0) at (0,3) {$\younganddescent{:14,:26,35}{124} \phimapsto$};
\node (1) at (1,3) {$\younganddescent{:25,:36,14}{235} \phimapsto$};
\node (2) at (2,3) {$\younganddescent{:23,:46,15}{346} \phimapsto$};
\node (3) at (3,3) {$\younganddescent{:14,:35,26}{145} \phimapsto$};
\node (4) at (4,3) {$\younganddescent{:12,:45,36}{256} \phimapsto$};
\node (5) at (4.9,3) {$\younganddescent{:13,:26,45}{136}$};
\draw [->, >=stealth] (4.9,3.42) -- (4.9,3.52) -- node[above]{\scriptsize $\phi$} (-0.1,3.52) -- (-0.1,3.42);

\node (6) at (0,2) {$\younganddescent{:15,:26,34}{125} \phimapsto$};
\node (7) at (1,2) {$\younganddescent{:12,:36,45}{236} \phimapsto$};
\node (8) at (2,2) {$\younganddescent{:13,:46,25}{134} \phimapsto$};
\node (9) at (3,2) {$\younganddescent{:24,:35,16}{245} \phimapsto$};
\node (10) at (4,2) {$\younganddescent{:23,:45,16}{356} \phimapsto$};
\node (11) at (4.9,2) {$\younganddescent{:14,:36,25}{146}$};
\draw [->, >=stealth] (4.9,2.42) -- (4.9,2.52) -- node[above]{\scriptsize $\phi$} (-0.1,2.52) -- (-0.1,2.42);

\node (12) at (0,1) {$\younganddescent{:13,:25,46}{135} \phimapsto$};
\node (13) at (1,1) {$\younganddescent{:12,:34,56}{246} \phimapsto$};
\node (14) at (2,1) {$\younganddescent{:13,:45,26}{135} \phimapsto$};
\node (15) at (3,1) {$\younganddescent{:24,:36,13}{246} \phimapsto$};
\node (16) at (4,1) {$\younganddescent{:15,:36,24}{135} \phimapsto$};
\node (17) at (4.9,1) {$\younganddescent{:12,:46,35}{246}$};
\draw [->, >=stealth] (4.9,1.42) -- (4.9,1.52) -- node[above]{\scriptsize $\phi$} (-0.1,1.52) -- (-0.1,1.42);

\node (18) at (1.5,0) {$\younganddescent{:14,:25,36}{1245} \phimapsto$};
\node (19) at (2.5,0) {$\younganddescent{:12,:35,46}{2356} \phimapsto$};
\node (20) at (3.4,0) {$\younganddescent{:13,:24,56}{1346}$};
\draw [->, >=stealth] (3.4,.42) -- (3.4,.52) -- node[above]{\scriptsize $\phi$} (1.4,.52) -- (1.4,.42);
\end{tikzpicture}
\end{ex}

Some remarks on this construction are in order.

\begin{rem}\label{rem:not-zn-action}
The bijection $\phi$ does not, in general, generate a $\ZZ_n$-action.
For example, the orbit of
\begin{equation*}
    T = \young(135\ten,248,679)
\end{equation*}
under the action of $\phi$ has order $20$.
\end{rem}

\begin{rem}
Because the operations $\RotSE, \pro, \RotNW$ commute with transposition,
$\phi$ also commutes with transposition, i.e.,
\begin{equation*}
    \phi(T^{\tra}) = (\phi T)^{\tra}
\end{equation*}
for all $T\in \SYT(\skewshape)$.
Moreover, note that
\begin{equation*}
    \phi^{-1} = \RotSE^{-1} \circ \dem \circ \RotNW.
\end{equation*}
Therefore, $\phi$ and $\phi^{-1}$ are conjugate with respect to reversal,
in the sense that
\begin{equation*}
    \phi^{-1} (T^{\rev}) = (\phi T)^{\rev}
\end{equation*}
for all $T\in \SYT(\skewshape)$.
\end{rem}

\section{Proof of Theorem~\ref{thm:main}}\label{sec:main-proof}

In this section we prove our main result.
Throughout this section, $\skewshape$ is a skew shape that is not a connected ribbon.
We will show that $(\cDes, \phi)$, defined in Theorem~\ref{thm:main},
has the three properties in Definition~\ref{defn:cyc-descent-map}.
The extension property holds by construction, so it remains to establish
the equivariance and non-Escher properties.

\subsection{Proving non-Escherness}

We first prove that $(\cDes, \phi)$ is non-Escher.
Note the following property of $\phi$.
\begin{lem}\label{lem:phi-increments-labels}
    Let $T\in \SYT(\skewshape)$,
    and let $\eta$ be some connected component of $\skewshape$.
    If $\eta$ has entries $J\subseteq [n]$ in $T$,
    then $\eta$ has entries $J+1$ in $\phi T$, where we take indices modulo $n$.
\end{lem}
\begin{proof}
    The operators $\RotSE$ and $\RotNW$ (and therefore, their inverses)
    only permute the cells within each connected component.
    The operator $\pro$ increments the entries of each connected component by $1$ modulo $n$.
\end{proof}

We will also use the following two-sided bound on $|\Des (T)|$.
\begin{lem}\cite[Lemma 4.6]{AER}\label{lem:des-size-bounds}
    Let $\skewshape$ be a skew shape with $n$ cells (possibly a connected ribbon).
    Let $c$ (resp. $r$) be the maximum length of a column (resp. row) in
    $\skewshape$.
    The maximum and minimum values of $|\Des(T)|$,
    over all $T\in \SYT(\skewshape)$, are
    \begin{equation*}
        \min |\Des(T)| = c-1,
        \quad
        \max |\Des(T)| = (n-1) - (r-1).
    \end{equation*}
\end{lem}

\begin{prop}\label{prop:non-escher}
    $(\cDes, \phi)$ satisfies the non-Escher property.
\end{prop}
\begin{proof}
    We will show $1\le |\cDes(T)|\le n-1$ for all $T\in \SYT(\skewshape)$.

    If any column of $\skewshape$ has height at least $2$, then by
    Lemma~\ref{lem:des-size-bounds}, $\min |\Des(T)| \ge 1$, so
    $\min |\cDes(T)| \ge 1$.
    Otherwise, $\skewshape$ is a disconnected union of single-row rectangles.
    Since $\skewshape$ is not a connected ribbon, there are at least
    two connected components.

    If $T\in \SYT(\skewshape)$ and $|\Des(T)|\ge 1$,
    then $|\cDes(T)|\ge 1$, as desired.
    Otherwise, $|\Des(T)|=0$.  Then, for each $i\in [n-1]$,
    $\pos_T(i)$ is (nonstrictly) south of $\pos_T(i+1)$.
    This implies that $\pos_T(1)$ is in the southernmost connected
    component and $\pos_T(n)$ is in the northernmost.
    By Lemma~\ref{lem:phi-increments-labels}, $\pos_{\phi T}(1)$ is
    in the northernmost connected component and $\pos_{\phi T}(2)$ is
    in the southernmost.
    Thus, $1\in \Des(\phi T)$, which implies $n\in \cDes(T)$.
    Hence, we always have $|\cDes(T)|\ge 1$.

    If any row of $\skewshape$ has length at least $2$, then by
    Lemma~\ref{lem:des-size-bounds}, $\max |\Des(T)| \le n-2$,
    so $\max |\cDes(T)| \le n-1$.
    Otherwise, $\skewshape$ is a disconnected union of single-column rectangles.
    Since $\skewshape$ is not a connected ribbon, there are at least
    two connected components.

    If $T\in \SYT(\skewshape)$ and $|\Des(T)|\le n-2$, then
    $|\cDes(T)|\le n-1$, as desired.
    Otherwise, $|\Des(T)|=n-1$.  So, for each $i\in [n-1]$,
    $\pos_T(i)$ is strictly north of $\pos_T(i+1)$.
    Thus, $\pos_T(1)$ is in the northernmost connected component,
    and $\pos_T(n)$ is in the southernmost.
    By Lemma~\ref{lem:phi-increments-labels}, $\pos_{\phi T}(1)$
    is in the southernmost connected component
    and $\pos_{\phi T}(2)$ is in the northernmost.
    Therefore $1\not\in \Des(\phi T)$ and $n\not\in \cDes(T)$.
    So, we always have $|\cDes(T)|\le n-1$.
\end{proof}

\subsection{Promotion Paths}
Before we prove equivariance, we will introduce some machinery about
promotion paths.

Recall that the \vocab{promotion path} of $T\in \SYT(\skewshape)$ is the set of cells
traversed by $1$ when $T$ is promoted.
\begin{defn}
    For $T\in \SYT(\skewshape)$, the \vocab{$\phi$-promotion path} of $T$,
    denoted $P_{\phi}(T)$, is the path traversed by $1$
    in the promotion phase of the action of $\phi$ on $T$.
    Equivalently, $P_{\phi}(T)$ is the promotion path of $\RotSE(T)$.
\end{defn}

\begin{ex}
    Consider the $T$ given in Example~\ref{ex:phi}.
    The promotion path $P_{\phi}(T)$ is shown in {\bf bold}.
    \begin{equation*}
        \young(:24,:35,16)
        \quad
        \overset{\RotSE}{\mapsto}
        \quad
        \young(:\btwo\bfour,:3\bsix,15)
        \quad
        \overset{\pro}{\mapsto}
        \quad
        \young(:13,:45,26)
        \quad
        \overset{\RotNW^{-1}}{\mapsto}
        \quad
        \young(:23,:45,16)
        .
    \end{equation*}
\end{ex}

This path starts at a southeast exterior corner of $\skewshape$,
which we call the path's \vocab{source},
and ends at a northwest exterior corner of $\skewshape$,
which we call the path's \vocab{destination}.
Note that the southeast and northwest exterior corners of $\skewshape$
both have a southwest-to-northeast order.

Throughout this subsection, we let $T\in \SYT(\skewshape)$ and
\begin{equation*}
    \RcSE(T) = \{n,\ldots,n-k+1\}.
\end{equation*}
Moreover, we let $\RcSE(T)$ appear in the southwest-to-northeast order
$x_1,\ldots,x_k$, such that
\begin{equation*}
    x_1>x_2>\cdots>x_i=n-k+1<x_{i+1}<\cdots<x_k.
\end{equation*}
We let $X$ be the southeast rotation endpoint of $T$ and $Y=\pos_T(n-k+1)=\pos_T(x_i)$.
We also define $T' = \pro(\RotSE(T))$.

Throughout this subsection, we take the convention that $\pos_T(x_0)$ and $\pos_T(x_{k+1})$ are,
respectively, the southwest and northeast corners of the connected component of $Y$.
In the following lemma, this convention is relevant when $i=1$ or $i=k$.
\begin{lem}\label{lem:drop-in-avoid}
    Suppose $n-k+1\in \RcSE(T')$, so $W = \pos_{T'}(n-k+1)$ is on the southeast boundary.
    Suppose further that $W$ is in the same connected component as $Y$.
    Then, $W$ is not on the portion of the southeast boundary
    between $\pos_T(x_{i-1})$ and $\pos_T(x_{i+1})$, inclusive.
\end{lem}
\begin{rem}
    If $k=1$, the portion of the southeast boundary between $\pos_T(x_0)$ and $\pos_T(x_{k+1})$
    is the portion of the southeast boundary in the connected component of $Y$.
    As a consequence, if the hypotheses of Lemma~\ref{lem:drop-in-avoid} hold, then $k>1$.
\end{rem}

\begin{proof}
    Because $1\in \RcNW(T)$, Lemma~\ref{lem:rot-no-interference} implies $1\not\in \RcSE(T)$.
    This implies $n-k+1>1$, so $n-k$ is a valid entry in $T$.
    Moreover, the action of $\RotSE$ on $T$ does not move $n-k$.
    This allows us to reason about $W' = \pos_T(n-k) = \pos_{\RotSE(T)}(n-k)$.
    Note that $W$ contains the incremented entry of $W'$ after promotion of $\RotSE(T)$.

    Note that $\RotSE$ rotates $n$ to $X$,
    and the promotion path $P_\phi(T)$ starts at $X$ and passes through $Y$.
    Thus in $T'$, the cells $\RpSE(T)\setminus \{Y\}$ collectively contain $\{n-k+2,\ldots,n\}$,
    the incremented entries of $\RcSE(T) \setminus \{n\}$.
    Since $W$ contains entry $n-k+1$ in $T$, this implies $W\not\in \RpSE(T)\setminus \{Y\}$.

    Let $Z_1 = \pos_{T}(x_{i-1})$ and $Z_2 = \pos_{T}(x_{i+1})$.
    We will show that $W$ is not northeast of $Z_1$ and southwest of $Y$.
    That $W$ is not northeast of $Y$ and southwest of $Z_2$ follows by an analogous argument.

    We handle the three possible configurations of $X$ and $Y$
    from Proposition~\ref{prop:rcse-locations} separately.

    Case 1: $X,Y$ are in configuration (a).

    Here, $Z_1$ is the southern neighbor of $Y$.
    Since $W\not\in \RpSE(T)\setminus \{Y\}$, $W\neq Z_1$. It suffices to show $W\neq Y$.

    Suppose $W=Y$; because $Y$ is on the promotion path $P_\phi(T)$,
    when $\RotSE(T)$ is promoted the incremented entry of $W'$ moves one cell east or south into $W$.
    Thus, $W'$ is the northern or western neighbor of $Y$.
    If $W'$ is the western neighbor of $Y$,
    the southern neighbor of $W'$ is not in $\RpSE(T)$ because $\RpSEc(T)$ is a hook;
    thus, the entry in this cell in $T$ is larger than $n-k$ but not in $\RcSE(T)$, contradiction.
    A similar contradiction arises if $W'$ is the northern neighbor of $Y$.

    Case 2: $X,Y$ are in configuration (b), with $X\neq Y$.

    Here, $Z_1$ is the southern neighbor of $Y$.
    Since $W\not\in \RpSE(T)\setminus \{Y\}$, $W\neq Z_1$. It again suffices to show $W\neq Y$.

    Suppose $W=Y$; as before, $W'$ is the northern or western neighbor of $Y$.
    If $W'$ is the western neighbor of $Y$,
    we get the same contradiction as the previous case by considering the southern neighbor of $W'$ in $T$.
    If $W'$ is the northern neighbor of $Y$, $W'$ is also on the southeast boundary,
    because otherwise $Y$ must be a southeast interior corner, and $X,Y$ would be in configuration (a).
    Thus, $\RcSE(T)\cup \{n-k\}$ is southeast min-unimodal in $T$,
    contradicting maximality of $k$.

    Case 3: $X,Y$ are in configuration (c), including possibly $X=Y$.

    Here, $Z_1$ is the northernmost cell of $\RpSEsw(T)$ if $\RpSEsw(T)$ is nonempty,
    and otherwise is the southwest corner of the connected component of $Y$.
    Suppose for contradiction that $W$ is on the southeast boundary, northeast of $Z_1$ and southwest of $Y$.

    If $W$ is strictly south of $Y$, it is not on the promotion path $P_{\phi}(T)$, so $W = W'$.
    Then, $\RcSE(T)\cup \{n-k\}$ is southeast min-unimodal in $T$, contradiction.

    Otherwise, $W$ is in the same row as $Y$.
    Recall that in $T'$, the cells $\RpSE(T)\setminus \{Y\}$ collectively contain $\{n-k+2,\ldots,n\}$,
    and $W$ contains $n-k+1$, the largest entry outside this set.
    Thus $W=Y$; as before, $W'$ is the northern or western neighbor of $Y$.
    If $W'$ is the northern neighbor of $Y$, we consider two sub-cases:
    \begin{enumerate}[label=(3\alph*)]
        \item If $X\neq Y$, we get the same contradiction as in Case 1 by considering the eastern neighbor of $W'$.
        \item If $X=Y$, $Y$ is a southeast exterior corner, so $W'$ is on the southeast boundary.
        Thus, $\RcSE(T)\cup \{n-k\}$ is southeast min-unimodal in $T$, contradicting maximality of $k$.
    \end{enumerate}
    If $W'$ is the western neighbor of $Y$, it is also on the southeast boundary,
    because otherwise $Y$ must be a southeast interior corner, and $X,Y$ would be in configuration (a).
    Thus, $\RcSE(T)\cup \{n-k\}$ is southeast min-unimodal in $T$, contradicting maximality of $k$.
\end{proof}

The main result of this subsection is the following.
\begin{prop}\label{prop:pro-source-location}
    Suppose $n-1$ and $n$ are in the same connected component of $T$, and $n-1\in \Des(T)$.
    Let $X_1, X_2$ denote the sources of $P_{\phi}(T)$ and $P_{\phi} (\phi T)$, respectively.
    If the first step in $P_{\phi}(T)$ is northward, $X_2$ is nonstrictly southwest of $X_1$.
    If the first step in $P_{\phi}(T)$ is westward, $X_2$ is strictly southwest of $X_1$.
\end{prop}

\begin{ex}
    Suppose
    \begin{equation*}
        T = \young(:256,139\ten,47\eleven,8\twelve).
    \end{equation*}
    The actions of $\phi$ on $T$ and $\phi T$ are shown below.
    The promotion paths $P_\phi(T)$ and $P_\phi(\phi T)$ are shown in {\bf bold}.
    \begin{eqnarray*}
        \young(:256,139\ten,47\eleven,8\twelve)
        \quad
        \overset{\RotSE}{\mapsto}
        \quad
        \young(:\btwo\bfive6,13\bnine\ten,47\btwelve,8\eleven)
        \quad
        \overset{\pro}{\mapsto}
        \quad
        \young(:137,246\eleven,58\ten,9\twelve)
        \quad
        \overset{\RotNW^{-1}}{\mapsto}
        \quad
        \young(:237,146\eleven,58\ten,9\twelve)
        &=& \phi T, \\
        \young(:237,146\eleven,58\ten,9\twelve)
        \quad
        \overset{\RotSE}{\mapsto}
        \quad
        \young(:237,\bone46\eleven,\bfive\beight\btwelve,9\ten)
        \quad
        \overset{\pro}{\mapsto}
        \quad
        \young(:348,157\twelve,269,\ten\eleven)
        \quad
        \overset{\RotNW^{-1}}{\mapsto}
        \quad
        \young(:348,157\twelve,269,\ten\eleven)
        &=& \phi^2 T.
    \end{eqnarray*}
    Here, the first step of $P_\phi(T)$ is northward, and $X_2$ is nonstrictly southwest of $X_1$.
\end{ex}

\begin{ex}
    Suppose
    \begin{equation*}
        T = \young(:256,1379,48\eleven,\ten\twelve).
    \end{equation*}
    The actions of $\phi$ on $T$ and $\phi T$ are shown below.
    The promotion paths $P_\phi(T)$ and $P_\phi(\phi T)$ are shown in {\bf bold}.
    \begin{eqnarray*}
        \young(:256,1379,48\eleven,\ten\twelve)
        \quad
        \overset{\RotSE}{\mapsto}
        \quad
        \young(:256,\bone379,\bfour\beight\btwelve,\ten\eleven)
        \quad
        \overset{\pro}{\mapsto}
        \quad
        \young(:367,148\ten,259,\eleven\twelve)
        \quad
        \overset{\RotNW^{-1}}{\mapsto}
        \quad
        \young(:367,148\ten,259,\eleven\twelve)
        &=& \phi T, \\
        \young(:367,148\ten,259,\eleven\twelve)
        \quad
        \overset{\RotSE}{\mapsto}
        \quad
        \young(:367,\bone48\ten,\btwo59,\beleven\btwelve)
        \quad
        \overset{\pro}{\mapsto}
        \quad
        \young(:478,159\eleven,26\ten,3\twelve)
        \quad
        \overset{\RotNW^{-1}}{\mapsto}
        \quad
        \young(:478,159\eleven,26\ten,3\twelve)
        &=& \phi^2 T.
    \end{eqnarray*}
    Here, the first step of $P_\phi(T)$ is westward, and $X_2$ is strictly southwest of $X_1$.
\end{ex}

For the rest of this subsection,
suppose $n-1$ and $n$ are in the same connected component of $T$, and $n-1\in \Des(T)$.
Since $n-1$ and $n$ are in the same connected component, $k > 1$.
Since $n-1\in \Des(T)$, $\pos_T(n)$ is southwesternmost in $\RpSE(T)$, so $x_1=n$.
Let the entry of $X$ in $T$ be $x_j$ for some $j\in [k]$.

Note that $X=X_1$, because both are the southeast rotation endpoint of $T$.

In $\RotSE(T)$, $x_1=n$ is at $X$,
and the entries $x_2,\ldots,x_k$ occupy the remaining cells of $\RpSE(T)$,
in that order, from southwest to northeast.
By examining the three configurations in Proposition~\ref{prop:rcse-locations},
we see that in $T'$,
the entries $x_2+1,\ldots,x_k+1$ are still in $\RpSE(T)$ in that order.
Moreover, $\{x_2+1,\ldots,x_k+1\}=\{n,\ldots,n-k+2\}$.  Thus,
\begin{equation*}
    \{n,\ldots,n-k+2\} \subseteq \RcSE(T').
\end{equation*}

By Lemma~\ref{lem:rot-no-interference}, the sets $\RcSE(\phi T)$ and
$\RcNW(\phi T)$ are disjoint.
Therefore each element of $\RcSE(\phi T)$ is in the same cell
in $\phi T$ and in
\begin{equation*}
    \RotNW(\phi T) = \pro(\RotSE(T)) = T'.
\end{equation*}
Moreover, Lemma~\ref{lem:rot-no-interference} states that
\begin{equation*}
    \RcSE(\phi T) = \RcSE(T').
\end{equation*}
Thus, $X_2$, the southeast rotation endpoint of $\phi T$,
is also the southeast rotation endpoint of $T'$.
We will use this characterization of $X_2$ in the following proof.

\begin{proof}[Proof of Proposition~\ref{prop:pro-source-location}]
    Recall that $X$ is due east or due south of $Y$.  We consider two cases.

    Case 1: $X$ is strictly due south of $Y$.

    This case corresponds to configuration (a) of Proposition~\ref{prop:rcse-locations},
    or configuration (b) where $X\neq Y$.
    In $\RotSE(T)$, $X$ contains the entry $n$;
    so, the beginning of $P_{\phi}(T)$ is a due-north path from $X$ to $Y$.
    We will show that $X_2$ is nonstrictly southwest of $X_1$.

    Since $X$ is strictly due south of $Y$, $j<i$.
    The action of $\RotSE$ on $T$ moves exactly the elements $x_1,\ldots,x_j$.
    In particular, it does not move $x_{j+1},\ldots,x_i$.
    Furthermore, promotion of $\RotSE(T)$ increments these elements and moves them down by one cell.
    So, for $m\in \{j+1,\ldots,i\}$, $\pos_{T'}(x_m+1)$ is the southern neighbor of $\pos_T(x_m)$.

    If $\RcSE(T')=\{n,\ldots,n-k+2\}$,
    $X_2$ is the southeast exterior corner due east or due south of $\pos_{T'}(n-k+2)=\pos_{T'}(x_i+1)$.
    Since $\pos_{T'}(x_i+1)$ is the southern neighbor of $\pos_T(x_i)=Y$, we have $X_2=X$.
    This is nonstrictly southwest of $X_1=X$, as desired.

    Otherwise, $n-k+1\in \RcSE(T')$,
    so $\{n,\ldots,n-k+1\}$ is southeast min-unimodal in $T'$.
    This implies that $\pos_{T'}(n-k+1)$ is on the southeast boundary
    between $\pos_{T'}(x_{i-1}+1)$ and $\pos_{T'}(x_{i+1}+1)$.
    Here, if $i=k$, we let $\pos_{T'}(x_{i+1}+1)$ be
    the northeast corner of the connected component of $Y$,
    and if $i=2$, we let $\pos_{T'}(x_{i-1}+1) = \pos_{T'}(n+1)$ be
    the southwest corner of the connected component of $Y$.

    We claim that $\pos_{T'}(x_{i+1}+1) = \pos_T(x_{i+1})$.  If $i=k$, this holds by definition.
    If $i<k$, the action of $\RotSE$ on $T$ does not move $x_{i+1}$.
    Moreover, in both configurations (a) and (b), $\pos_T(x_{i+1})$ is strictly east of $X$;
    therefore $\pos_T(x_{i+1})$ is not in $P_{\phi}(T)$,
    and promotion of $\RotSE(T)$ increments its entry without moving it.
    Therefore $\pos_{T'}(x_{i+1}+1) = \pos_T(x_{i+1})$.

    If $i>j+1$, then $i-1\in \{j+1,\ldots,i\}$.
    So, $\pos_{T'}(x_{i-1}+1)$ is the southern neighbor of $\pos_T(x_{i-1})$.
    Therefore, $\pos_{T'}(n-k+1)$ is on the southeast boundary
    between the southern neighbor of $\pos_T(x_{i-1})$ and $\pos_T(x_{i+1})$.
    But, by Lemma~\ref{lem:drop-in-avoid},
    $\pos_{T'}(n-k+1)$ is not on the southeast boundary between $\pos_T(x_{i-1})$ and $\pos_T(x_{i+1})$.
    Thus, $\pos_{T'}(n-k+1)$ is the southern neighbor of $\pos_T(x_{i-1})$, which is $\pos_{T'}(x_{i-1}+1)$.
    But $x_{i-1}+1 > n-k+1$, contradiction.

    Otherwise, $i=j+1$.
    Define $Z$ to be the northernmost cell of $\RpSEsw(T)$ if it is nonempty,
    and otherwise as the southwest corner of the connected component of $Y$.
    We claim that $\pos_{T'}(x_{i-1}+1) = Z$.
    If $j>1$, $\RpSEsw(T)$ is nonempty, so $Z$ is the northernmost cell of $\RpSEsw(T)$.
    The action of $\RotSE$ moves $x_j=x_{i-1}$ to $Z$.
    Since $Z$ is strictly south of $X$, it is not on $P_\phi(T)$;
    thus, promotion of $\RotSE(T)$ increments the entry in $Z$ without moving it, and $\pos_{T'}(x_{i-1}+1) = Z$.
    If $j=1$, $\RpSEsw(T)$ is empty, so $Z$ is the southwest corner of the connected component of $Y$.
    As $\pos_{T'}(x_{i-1}+1)$ equals this cell by definition, $\pos_{T'}(x_{i-1}+1) = Z$.

    So, $\pos_{T'}(n-k+1)$ is on the southeast boundary between $Z$ and $\pos_T(x_{i+1})$.
    By Lemma~\ref{lem:drop-in-avoid},
    $\pos_{T'}(n-k+1)$ is not on the southeast boundary
    between $\pos_T(x_{i-1})=X$ and $\pos_T(x_{i+1})$.
    Thus, $\pos_{T'}(n-k+1)$ is southwest of $X$.
    By southeast min-unimodality,
    $\pos_{T'}(\min(\RcSE(T')))$ is southwest of $X$ as well.
    Hence, $X_2$ is nonstrictly southwest of $X_1=X$.

    Case 2: $X$ is due east of $Y$, including possibly $X=Y$.

    This case corresponds to configuration (c) of Proposition~\ref{prop:rcse-locations}.
    We will show that in this case, $X_2$ is strictly southwest of $X_1$.
    Note that in this case the first step of $P_\phi(T)$ may be northward or westward,
    and in particular may be northward when $X=Y$;
    but, regardless of the direction of this step,
    proving that $X_2$ is strictly southwest of $X_1$ implies the result.

    We claim that in this case, $\pos_T(n) \in \RpSEsw(T)$.
    Suppose otherwise; since $\pos_T(n)$ is southwesternmost in $\RpSE(T)$, $\pos_T(n)\in \RpSEc(T)$.
    Since $n$ is the largest entry, $\pos_T(n) = X$.
    But, since $\pos_T(n)$ is southwesternmost in $\RpSE$, $X=Y$.
    This is only possible if $n=n-k+1$, which implies $k=1$.
    This is a contradiction because $n-1$ and $n$ are in the same connected component of $T$,
    which implies $k>1$.
    As a consequence, $\RpSEsw(T)$ is nonempty and $i>1$.
    Let $Z$ be the northernmost cell of $\RpSEsw(T)$.

    In $\RotSE(T)$, $X$ contains $n$ and $Z$ contains $n-k+1$.
    Since $Z$ is strictly southwest of $X$, it does not intersect $P_{\phi}(T)$.
    Thus, promoting $\RotSE(T)$ increments but does not move the entry in $Z$,
    and $\pos_{T'}(n-k+2)=Z$.

    If $\RcSE(T')=\{n,\ldots,n-k+2\}$, then $X_2$ is
    a southeast exterior corner due east or due south of $\pos_{T'}(n-k+2)=Z$.
    As $Z$ is strictly southwest of $X$, $X_2$ is strictly southwest of $X_1=X$.

    Otherwise, $n-k+1\in \RcSE(T')$,
    so $\{n,\ldots,n-k+1\}$ is southeast min-unimodal in $T'$.
    This implies that $\pos_{T'}(n-k+1)$ is on the southeast boundary,
    between $\pos_{T'}(x_{i-1}+1)$ and $\pos_{T'}(x_{i+1}+1)$.
    As before, if $i=k$, we let $\pos_{T'}(x_{i+1}+1)$ be the northeast corner
    of the connected component of $Y$; as noted above, we never have $i=1$.

    We claim that $\pos_{T'}(x_{i+1}+1) = \pos_T(x_{i+1})$.
    If $i=k$, this is true by definition; otherwise, assume $i<k$.
    If $\pos_T(x_{i+1}) \in \RpSEc(T)$, $\RotSE$ moves the entry in this cell one cell westward;
    thus, $\pos_{\RotSE(T)}(x_{i+1})$ is the western neighbor of $\pos_T(x_{i+1})$.
    This cell lies on the promotion path of $\RotSE(T)$,
    so $\pro$ increments this cell and moves its incremented entry one cell eastward.
    So, $\pos_{T'}(x_{i+1}+1) = \pos_T(x_{i+1})$.
    If $\pos_T(x_{i+1}) \in \RpSEne(T)$, $\RotSE$ does not move this cell's entry,
    so $\pos_{\RotSE(T)}(x_{i+1}) = \pos_T(x_{i+1})$.
    This cell does not lie on the promotion path of $\RotSE(T)$,
    so $\pro$ increments this cell but does not move its incremented entry.
    Thus, again $\pos_{T'}(x_{i+1}+1) = \pos_T(x_{i+1})$.

    By Lemma~\ref{lem:drop-in-avoid}, $\pos_{T'}(n-k+1)$ is not on the southeast boundary between the cells
    $\pos_T(x_{i-1})=Z$ and $\pos_T(x_{i+1})=\pos_{T'}(x_{i+1}+1)$.
    Thus, $\pos_{T'}(n-k+1)$ is southwest of $Z$.
    By southeast min-unimodality, $X_2 = \pos_{T'}(\min(\RcSE(T')))$ is southwest of $Z$ as well.
    Since $Z$ is strictly southwest of $X=X_1$, $X_2$ is strictly southwest of $X_1$.
\end{proof}

\subsection{Demotion Paths}

We analogously define a notion of demotion paths.

The \vocab{demotion path} of $T\in \SYT(\skewshape)$ is the set of cells
traversed by $n$ when $T$ is demoted.
\begin{defn}
    For $T\in \SYT(\skewshape)$, the \vocab{$\phi^{-1}$-demotion path} of $T$,
    denoted $P_{\phi^{-1}}(T)$, is the path traversed by $n$
    in the demotion phase of the action of $\phi^{-1}$ on $T$.
    Equivalently, $P_{\phi^{-1}}(T)$ is the demotion path of $\RotNW(T)$.
\end{defn}

This path starts at a northwest exterior corner of $\skewshape$,
which we call the path's \vocab{source},
and ends at a southeast exterior corner,
which we call its \vocab{destination}.
Note that for all $T$, the $\phi^{-1}$-demotion path of $\phi T$
is the $\phi$-promotion path of $T$, traversed in the opposite direction.

We have this analogue of Proposition~\ref{prop:pro-source-location}.
\begin{cor}\label{cor:dem-source-location}
    Suppose $1$ and $2$ are in the same connected component of $T$, and $1\not\in \Des(T)$.
    Let $X_1, X_2$ denote the sources of $P_{\phi^{-1}}(T)$ and $P_{\phi^{-1}}(\phi^{-1} T)$, respectively.
    If the first step in $P_{\phi^{-1}}(T)$ is southward, $X_2$ is strictly southwest of $X_1$.
    If the first step in $P_{\phi^{-1}}(T)$ is eastward, $X_2$ is nonstrictly southwest of $X_1$.
\end{cor}
\begin{proof}
    Consider the tableau $T' = (T^{\rev})^{\tra}$.
    Since $1$ and $2$ are in the same connected component of $T$,
    $n-1$ and $n$ are in the same connected component of $T'$.
    Since $1\not\in \Des(T)$, $n-1\in \Des(T')$.

    The promotion paths $P_{\phi}(T')$ and $P_{\phi}(\phi T')$ are
    the reverse transposes of
    the demotion paths $P_{\phi^{-1}}(T)$ and $P_{\phi^{-1}}(\phi^{-1} T)$,
    respectively.
    Let $X_1'$ and $X_2'$ be the sources of $P_{\phi}(T')$ and $P_{\phi}(\phi T')$.
    $X_1'$ and $X_2'$ are the reverse transposes of $X_1$ and $X_2$.

    Suppose the first step in the $P_{\phi^{-1}}(T)$ is southward.
    Then, the first step in $P_{\phi}(T')$ is westward.
    By Proposition~\ref{prop:pro-source-location} on $T'$,
    $X_2'$ is strictly southwest of $X_1'$.
    Thus, $X_2$ is strictly southwest of $X_1$.

    Conversely, suppose the first step in the $P_{\phi^{-1}}(T)$ is eastward.
    Then, the first step in $P_{\phi}(T')$ is northward.
    By Proposition~\ref{prop:pro-source-location} on $T'$,
    $X_2'$ is nonstrictly southwest of $X_1'$.
    Thus, $X_2$ is nonstrictly southwest of $X_1$.
\end{proof}

\subsection{Pseudo-Promotion Paths}

We will also use a notion of pseudo-promotion paths.
\begin{defn}
    For $T\in \SYT(\skewshape)$ and a southeast exterior corner $Z$ of $\skewshape$,
    the \vocab{pseudo-promotion path} of $Z$ in $T$ is the set of cells starting at $Z$,
    obtained by the following procedure:
    while the current cell is not a northwest exterior corner,
    walk to the larger of the current cell's northern neighbor (if it exists)
    and western neighbor (if it exists).
\end{defn}
\begin{ex}\label{ex:pseudo-promotion}
    For the tableau $T$ below, the pseudo-promotion path
    for the cell containing $11$ is shown in {\bf bold}.
    \begin{equation*}
        \young(:167,\btwo\bfive89,3\bten\beleven,4\twelve).
    \end{equation*}
\end{ex}

If $\pos_T(n)=Z$, the pseudo-promotion path of $Z$ in $T$ is the promotion path of $T$.
However, the pseudo-promotion path is defined even when the entry in $Z$ is not $n$.

\begin{lem}\label{lem:rotse-same-pseudo-promotion-path}
    Let $T\in \SYT(\skewshape)$ and $Z$ be a southeast exterior corner of $\skewshape$.
    The pseudo-promotion paths of $Z$ in $T$ and $\RotSE(T)$ are the same.
\end{lem}
\begin{proof}
    For each cell $W$ that is not a northwest exterior corner, define the
    \vocab{pseudo-promotion direction} of $W$ in $T$
    as the direction of the larger of $W$'s northern neighbor (if it exists)
    and its western neighbor (if it exists).
    This is the direction a pseudo-promotion path would take to leave $W$.

    We claim that for any $W$,
    the pseudo-promotion directions of $W$ in $T$ and $\RotSE(T)$ are the same.

    In both $T$ and $\RotSE(T)$, the entry in $\RpSE(T)$ is $\RcSE(T)$.
    Proposition~\ref{prop:rcse-locations} implies that $\RpSE(T)$ avoids the shape
    \begin{equation*}
        \young(:\blank,\blank\blank).
    \end{equation*}
    Thus, it is not possible for both the northern and western neighbors of $W$
    to be in $\RpSE(T)$.

    We consider two cases.

    Case 1: One of the northern and western neighbors of $W$ is in $\RpSE(T)$.

    Since all elements of $\RcSE(T)$ are larger than all non-elements of
    $\RcSE(T)$, the pseudo-promotion direction of $W$ in both $T$ and
    $\RotSE(T)$ is in the direction of the neighbor in $\RpSE(T)$.

    Case 2: Neither of the northern and western neighbors of $W$ is in $\RpSE(T)$.

    The operation $\RotSE$ does not move either of these neighbors, so the northern
    and western neighbors of $W$ have the same entries in $T$ and $\RotSE(T)$.
    Therefore the pseudo-promotion direction of $W$ is the same in both.
\end{proof}

If a pseudo-promotion path reaches a northwest interior corner $U$ of $\skewshape$,
the destination of the path must be the northwest exterior corner due north of $U$
or the northwest exterior corner due west of $U$.
So, the remaining steps of the path must be all northward or all westward.

\begin{lem}\label{lem:rotnw-almost-same-pseudo-promotion-path}
    Let $T\in \SYT(\skewshape)$, and let $Z$ be a southeast exterior corner of $\skewshape$.
    If the pseudo-promotion path of $Z$ in $T$ passes through a northwest interior corner $U$,
    the pseudo-promotion paths of $Z$ in $T$ and $\RotNW(T)$ are equal between $Z$ and $U$, inclusive.
    Otherwise, the pseudo-promotion paths of $Z$ in $T$ and $\RotNW(T)$ are equal.
\end{lem}
Thus, the two paths are not always equal;
but if they differ, they differ only in that after reaching $U$,
one travels due north while the other travels due west.
\begin{ex}\label{ex:pseudo-promotion-rotnw}
    Let $T$ be the tableau in Example~\ref{ex:pseudo-promotion}.
    The tableau $\RotNW(T)$ is shown below, with the pseudo-promotion path
    for the cell containing $11$ shown in {\bf bold}.
    \begin{equation*}
        \young(:\btwo67,1\bfive89,3\bten\beleven,4\twelve).
    \end{equation*}
    This path differs from the pseudo-promotion path in Example~\ref{ex:pseudo-promotion}
    only in that after reaching the northwest interior corner with entry $5$,
    this path traveled north while the path in Example~\ref{ex:pseudo-promotion} traveled west.
\end{ex}

\begin{proof}
    We will show that for any cell $W$ not a northwest interior corner,
    the pseudo-promotion directions of $W$ in $T$ and $\RotNW(T)$ are the same.

    We claim that not both the northern and western neighbors of $W$ are in $\RpNW(T)$.
    Suppose otherwise; because $W$ has northern and western neighbors
    and is not a northwest interior corner, it has a northwestern neighbor.
    As the northern and western neighbors of $W$ are in $\RpNW(T)$,
    the northwestern neighbor is also in $\RpNW(T)$.
    Thus, $\RpNW(T)$ contains the shape
    \begin{equation*}
        \young(\blank\blank,\blank:).
    \end{equation*}
    But, Proposition~\ref{prop:rcnw-locations} implies that $\RpNW(T)$ avoids this shape,
    contradiction.

    So, there are two cases to consider.

    Case 1: One of the northern and western neighbors of $W$ is in $\RpNW(T)$.

    Since all elements of $\RcNW(T)$ are smaller than all non-elements of
    $\RcNW(T)$, the pseudo-promotion direction of $W$ in both $T$ and $\RotNW(T)$ is the same.

    Case 2: Neither of the northern and western neighbors of $W$ is in $\RpNW(T)$.

    The operation $\RotNW$ does not move either of these neighbors,
    so the pseudo-promotion direction of $W$ is the same in both $T$ and $\RotNW(T)$.
\end{proof}

\subsection{Cyclic Rotation of Descents}

The main result of this subsection is the following.
\begin{prop}\label{prop:double-promotion}
    Let $T\in \SYT(\skewshape)$.  If $n-1 \in \Des(T)$, then $1\in \Des(\phi^2 T)$.
\end{prop}

Throughout this subsection, let
\begin{equation*}
    P_1 = P_{\phi}(T), \quad P_2 = P_{\phi}(\phi T).
\end{equation*}
Let $S_1, S_2$ denote the sources of $P_1$ and $P_2$,
and let $D_1, D_2$ denote their destinations, respectively.

By definition, $P_2$ is the pseudo-promotion path of $S_2$ in $\RotSE(\phi T)$.
By Lemma~\ref{lem:rotse-same-pseudo-promotion-path} applied
to $\phi T$, $P_2$ is also the pseudo-promotion path of $S_2$ in $\phi T$.

Let
\begin{equation*}
	T' = \RotNW(\phi T) = \pro(\RotSE(T)).
\end{equation*}
Let $P_2'$ be the pseudo-promotion path of $S_2$ in $T'$, with destination $D_2'$.

\begin{prop}\label{prop:pro-dest-location}
    Suppose $1$ and $2$ are in the same connected component of $\phi^2 T$, and $1\not \in \Des(\phi^2 T)$.
    If the last step in $P_2'$ is northward, then $D_1$ is strictly southwest of $D_2'$.
    If the last step in $P_2'$ is westward, then $D_1$ is nonstrictly southwest of $D_2'$.
\end{prop}
\begin{ex}
    Suppose
    \begin{equation*}
        T = \young(:138,247\twelve,59\eleven,6\ten).
    \end{equation*}
    The action of $\phi$ on $T$ is shown below.
    The promotion path $P_1$ and pseudo-promotion path $P_2'$ are shown in {\bf bold} on,
    respectively, $\RotSE(T)$ and $T' = \pro(\RotSE(T))$.
    Note that $P_2'$ begins on the cell with entry $10$ in $T'$
    because this is the southeast rotation endpoint of $\phi T$.
    \begin{eqnarray*}
        \young(:138,247\twelve,59\eleven,6\ten)
        \quad
        \overset{\RotSE}{\mapsto}
        \quad
        \young(:138,\btwo47\eleven,\bfive\bnine\btwelve,6\ten)
        \quad
        \overset{\pro}{\mapsto}
        \quad
        \young(:\btwo49,1\bfive\beight\twelve,36\bten,7\eleven)
        \quad
        \overset{\RotNW^{-1}}{\mapsto}
        \quad
        \young(:149,258\twelve,36\ten,7\eleven).
    \end{eqnarray*}
    The last step in $P_2'$ is northward, and $D_1$ is strictly southwest of $D_2'$.
\end{ex}
\begin{ex}
    Suppose
    \begin{equation*}
        T = \young(:14\ten,238\twelve,57\eleven,69).
    \end{equation*}
    The action of $\phi$ on $T$ is shown below.
    The promotion path $P_1$ and pseudo-promotion path $P_2'$ are shown in {\bf bold} on,
    respectively, $\RotSE(T)$ and $T' = \pro(\RotSE(T))$.
    Note that $P_2'$ begins on the cell with entry $12$ in $T'$
    because this is the southeast rotation endpoint of $\phi T$.
    \begin{eqnarray*}
        \young(:14\ten,238\twelve,57\eleven,69)
        \quad
        \overset{\RotSE}{\mapsto}
        \quad
        \young(:\bone\bfour\ten,23\beight\eleven,57\btwelve,69)
        \quad
        \overset{\pro}{\mapsto}
        \quad
        \young(:\bone\btwo\beleven,345\btwelve,689,7\ten)
        \quad
        \overset{\RotNW^{-1}}{\mapsto}
        \quad
        \young(:12\eleven,345\twelve,689,7\ten).
    \end{eqnarray*}
    The last step in $P_2'$ is westward, and $D_1$ is nonstrictly southwest of $D_2'$.
\end{ex}

\begin{proof}
    Recall that $P_{\phi^{-1}}(\phi T)$ and $P_{\phi^{-1}}(\phi^2 T)$ are,
    respectively, $P_1$ and $P_2$ traversed in the opposite direction.
    Corollary~\ref{cor:dem-source-location}, applied to $\phi^2 T$, implies that
    if the last step in $P_2$ is northward,
    $D_1$ is strictly southwest of $D_2$,
    and if the last step in $P_2$ is westward,
    $D_1$ is nonstrictly southwest of $D_2$.

    If $P_2'=P_2$, then $D_2=D_2'$ and the result follows.

    Otherwise, by Lemma~\ref{lem:rotnw-almost-same-pseudo-promotion-path},
    $P_2'$ and $P_2$ are the same until they reach a northwest interior
    corner $C$, and then one of them moves north while the other moves west.

    If $P_2$ moves north and $P_2'$ moves west after $C$, then $D_2'$
    is the first northwest exterior corner to the southwest of $D_2$.
    Since $D_1$ is strictly southwest of $D_2$,
    $D_1$ is nonstrictly southwest of $D_2'$.

    If $P_2$ moves west and $P_2'$ moves north after $C$, then $D_2'$
    is the first northwest exterior corner to the northeast of $D_2$.
    Since $D_1$ is nonstrictly southwest of $D_2$,
    $D_1$ is strictly southwest of $D_2'$.
\end{proof}

\begin{proof}[Proof of Proposition~\ref{prop:double-promotion}]
    The operations $\RotSE$ and $\RotNW$ do not move entries across connected components.
    The operation $\pro$ increments all entries, but does not move any incremented entries
    across connected components.
    Thus, the connected components of $n-1$ and $n$ in $T$ are, respectively,
    the connected components of $1$ and $2$ in $\phi^2 T$.

    Let us first address when $n-1$ and $n$ are not in the same connected component of $T$.
    If the connected component of $n-1$ is north of the connected component of $n$ in $T$,
    the connected component of $1$ is north of the connected component of $2$ in $\phi^2 T$.
    Consequently, $n-1\in \Des(T)$ and $1\in \Des(\phi^2T)$.
    Analogously, if the connected component of $n-1$ is south of the connected component of $n$ in $T$,
    then $n-1\not\in \Des(T)$ and $1\not\in \Des(\phi^2 T)$.
    This proves the result when $n-1$ and $n$ are not in the same connected component of $T$.

    Otherwise, assume $n-1$ and $n$ are in the same connected component of $T$.
    It follows that $1$ and $2$ are in the same connected component of $\phi^2 T$.
    Suppose for contradiction that $n-1\in \Des(T)$ and $1\not\in \Des(\phi^2 T)$.

    By Propositions~\ref{prop:pro-source-location} and
    \ref{prop:pro-dest-location}, we have the following conditions.
    \begin{enumerate}[label=(\alph*)]
        \item If the first step in $P_1$ is northward, then $S_2$ is nonstrictly southwest of $S_1$.
        \item If the first step in $P_1$ is westward, then $S_2$ is strictly southwest of $S_1$.
        \item If the last step in $P_2'$ is northward, then $D_1$ is strictly southwest of $D_2'$.
        \item If the last step in $P_2'$ is westward, then $D_1$ is nonstrictly southwest of $D_2'$.
    \end{enumerate}

    Let
    \begin{equation*}
        \Gamma = \{Z\in P_1 | \text{western neighbor of $Z$ is also in $P_1$}\}
        \cup \{D_1\}.
    \end{equation*}
    Each column that $P_1$ intersects contains one cell in $\Gamma$.

    We claim that if $Z\in \Gamma \cap P_2'$, then the southern neighbor of $Z$ is not in $P_2'$.
    Suppose for contradiction that $Z\in \Gamma\cap P_2'$ and the southern neighbor of $Z$ is in $P_2'$.
    Since $Z\in \Gamma$, $Z=D_1$ or the western neighbor of $Z$ is in $P_1$.

    Suppose $Z=D_1$.
    Then, the southern neighbor of $D_1$ is in $P_2'$, so $D_2'$ is nonstrictly southwest of $D_1$.
    But, by conditions (c) and (d), $D_1$ is nonstrictly southwest of $D_2'$. Thus $D_1=D_2'$.
    Since the southern neighbor of $D_1=D_2'$ is in $P_2'$, the last step in $P_2'$ is northward.
    This contradicts condition (c), that $D_1$ is strictly southwest of $D_2'$.

    Otherwise, the western neighbor of $Z$ is in $P_1$. Because the southern neighbor of $Z$ exists,
    the southwestern neighbor of $Z$ must also exist, or else $\skewshape$ is not a skew shape.
    So, in $\RotSE(T)$, let the western neighbor and southwestern neighbors of $Z$
    have entries $u$ and $v$, as shown below.
    \begin{equation*}
        \young(uZ,v\blank).
    \end{equation*}
    Thus, $u<v$.
    Since both $Z$ and its western neighbor are in $P_1$, the entries in these cells in $T'$ are
    \begin{center}
        \begin{tabular}{|c|c|}
            \hline
            & u+1 \\ \hline
            v+1 & \\ \hline
        \end{tabular}.
    \end{center}
    Since both $Z$ and its southern neighbor are in $P_2'$,
    the step in $P_2'$ departing the southern neighbor of $Z$ is northward, contradicting that $u<v$.
    This proves our claim that the southern neighbor of $Z\in \Gamma\cap P_2'$ cannot be in $P_2'$.

    The path $P_2'$ runs southeast to northwest.
    Since $S_2$ is nonstrictly southwest of $S_1$ by conditions (a) and (b),
    and $D_2'$ is nonstrictly northeast of $D_1$ by conditions (c) and (d),
    $P_2'$ must be confined to the region of $\skewshape$ between the columns of $S_1$ and $D_1$, inclusive.
    Call this region $\Delta$.
    As each column of $\Delta$ contains one element of $\Gamma$, $\Gamma$ divides $\Delta$ into regions
    $\Delta^+$, the cells nonstrictly due north of a cell in $\Gamma$,
    and $\Delta^-$, the cells strictly due south of a cell in $\Gamma$.

    Conditions (a) and (b) imply $S_2\in \Delta^-$, while conditions (c) and (d) imply $D_2' \in \Delta^+$.
    Thus, $P_2'$ must intersect $\Gamma$ at a cell $Z$ such that the southern neighbor of $Z$ is also in $P_2'$.
    This is a contradiction.
\end{proof}

\subsection{Completion of the Proof}
We now have the tools to show that $(\cDes, \phi)$
satisfies the equivariance property.

\begin{prop}\label{prop:equivariance}
$(\cDes, \phi)$ satisfies the equivariance property.
That is, if $T\in \SYT(\skewshape)$ and $i\in[n]$, then $i\in \cDes(T)$
if and only if $i+1\in \cDes(\phi T)$, where we take indices modulo $n$.
\end{prop}
\begin{proof}
We consider three cases: $i\in [n-2]$, $i=n-1$, and $i=n$.

Case 1: $i\in [n-2]$.

The result follows from Proposition~\ref{prop:rotse-preserve-descents},
Lemma \ref{lem:pro-increments-des},
and Proposition \ref{prop:rotnw-preserve-descents},
due to the following equivalence chain:
\begin{align*}
    i\in \Des(T)
    & \Leftrightarrow
    i\in \Des(\RotSE(T)) \\
    & \Leftrightarrow
    i+1\in \Des(\pro(\RotSE(T))) \\
    & \Leftrightarrow
    i+1\in \Des(\RotNW^{-1}(\pro(\RotSE(T)))).
\end{align*}

Case 2: $i=n-1$.

It is equivalent to show that $n-1\in \Des(T)$ if and only if
$1\in \Des(\phi^2 T)$.

The forward direction follows from Proposition~\ref{prop:double-promotion}.

To show the converse, we consider $T$ such that $n-1\not\in \Des(T)$.
Then, $n-1\in \Des(T^{\tra})$.
By Proposition~\ref{prop:double-promotion}, $1\in \Des(\phi^2 (T^{\tra}))$.
Since $\phi$ commutes with transposition, $1\in \Des((\phi^2 T)^{\tra})$,
and thus, $1\not \in \Des(\phi^2 T)$.
Therefore, $n-1\not\in \Des(T)$ implies $1\not\in \Des(\phi^2 T)$, as desired.

Case 3: $i=n$.

By the definition of $\cDes$, $n\in \cDes(T)$
if and only if $1\in \cDes(\phi T)$.
\end{proof}

Finally, we prove our main result.
\begin{proof}[Proof of Theorem~\ref{thm:main}]
We will show that $(\cDes,\phi)$ satisfies the
extension, equivariance, and non-Escher properties.
Extension follows by definition.  Equivariance and non-Escher follow from
Propositions~\ref{prop:equivariance} and \ref{prop:non-escher}.
\end{proof}

\section{Special Cases}\label{sec:special-cases}

In this section we review some of the literature
for cyclic descents of skew tableaux.
We will show that several results in the literature
are special cases of Theorem~\ref{thm:main}.

First, let us introduce the relevant definitions from the literature.

For a (not necessarily standard) tableau $T$ of shape $\skewshape$ and integer $a$,
let $a+T$ be the tableau obtained by adding $a$ to each entry of $T$ and taking entries modulo $n$.

Where appropriate, we will also name cells by their coordinates.
Thus, $(x,y)$ is the cell in the $x^{\text{th}}$ row and $y^{\text{th}}$ column of $\skewshape$
where $(1,1)$ is the shared northwest corner of $\lambda$ and $\mu$.
Let $T_{x,y}$ denote the entry of $T$ in cell $(x,y)$.

The \vocab{generalized jeu de taquin} operator $\jdt$ \cite[Definition 2.3]{AER} acts on a tableau $T$ as follows.
While $T$ is not standard, apply the following operation to $T$:
find the minimal entry $i$ such that $i$ is not smaller than both its northern and western neighbors,
and swap $i$ with the larger of these neighbors.

\begin{ex}
    Below is one action of $\jdt$ on a nonstandard tableau.
    \begin{equation*}
        \young(342,51)
        \quad \mapsto \quad
        \young(342,15)
        \quad \mapsto \quad
        \young(142,35)
        \quad \mapsto \quad
        \young(124,35).
    \end{equation*}
\end{ex}
\begin{rem}
    For all $T\in \SYT(\skewshape)$,
    \begin{equation*}
    \pro(T) = \jdt(1+T)
    \quad
    \text{and}
    \quad
    \dem(T) = \jdt(-1+T).
    \end{equation*}
\end{rem}

Throughout this section, we adopt the following notational convention.
We use $(\cDes, \phi)$ to denote the cyclic descent map given in Theorem~\ref{thm:main},
and $(\cDes',\phi')$ to denote the cyclic descent maps from the literature
that we are comparing to $(\cDes, \phi)$.

The following lemma allows us to show that two cyclic descent maps
$(\cDes,\phi)$ and $(\cDes',\phi')$ coincide by
just showing that $\phi$ and $\phi'$ coincide.
\begin{lem}\label{lem:phi-sufficient}
    If a bijection $\phi:\SYT(\skewshape)\rightarrow \SYT(\skewshape)$
    is part of a cyclic descent map $(\cDes,\phi)$,
    then $\cDes$ is uniquely determined by $\phi$.
\end{lem}
\begin{proof}
    The action of $\phi$ partitions $\SYT(\skewshape)$ into orbits.
    Within each orbit, $|\cDes(T)|$ is constant.
    Thus, within each orbit $|\Des(T)|$ takes up to two distinct values:
    $|\cDes(T)|-1$ if $n\in \cDes(T)$, and otherwise $|\cDes(T)|$.
    In fact, by the non-Escher property,
    each orbit contains $T$ such that $n\in \cDes(T)$ and $T$ such that $n\not\in \cDes(T)$.
    So, within each orbit $|\Des(T)|$ takes exactly two distinct values.

    Therefore, the function $\cDes$ is uniquely determined from $\phi$ as follows:
    if, among the two descent set sizes in the orbit of $T$,
    $|\Des(T)|$ is the smaller value, let $n\in \cDes(T)$;
    if $|\Des(T)|$ is the larger value, let $n\not\in \cDes(T)$.
\end{proof}

\subsection{Rectangles and Strips}

Recall this result by Rhoades, stated in Section~\ref{sec:intro}.
\begin{repthm}{thm:rhoades}\cite[Lemma 3.3]{Rho}
    Let $\skewshape$ be a rectangular Young diagram with length and width both larger than 1.
    Let $\phi' = \pro$, and let $n\in \cDes'(T)$ if and only if $n-1 \in \Des(\dem T)$.
    Then $(\cDes', \phi')$ is a cyclic descent map.
\end{repthm}

Adin, Elizalde, and Roichman give an analogous construction for \vocab{strips},
skew shapes with more than one connected component
whose connected components are all single-row or single-column rectangles.
\begin{thm}\label{thm:aer-strip}\cite[Proposition 3.3]{AER}
    Let $\skewshape$ be a strip.
    Let $\phi' = \pro$, and let $n\in \cDes'(T)$ if and only if $n-1 \in \Des(\dem T)$.
    Then $(\cDes', \phi')$ is a cyclic descent map.
\end{thm}

\begin{prop}
    The cyclic descent map $(\cDes', \phi')$,
    defined in Theorems~\ref{thm:rhoades} and \ref{thm:aer-strip},
    coincides with the $(\cDes, \phi)$ defined in Theorem~\ref{thm:main}.
\end{prop}
\begin{proof}
    A rectangle has only one northwest exterior corner and one southeast exterior corner.
    Thus, in any standard Young tableau of a rectangle,
    the northwest and southeast rotation endpoints are always these corners.
    Therefore the actions of both $\RotNW$ and $\RotSE$ on rectangles are the identity, and
    \begin{equation*}
        \phi = \RotNW^{-1} \circ \pro \circ \RotSE = \pro = \phi'.
    \end{equation*}

    On any standard Young tableau of a strip,
    the actions of $\RotNW$ and $\RotSE$ are similarly the identity,
    because the northwest rotation corner is always the cell containing $1$
    and the southeast rotation corner is always the cell containing $n$.
    Thus, $\phi = \pro = \phi'$ on strips as well.

    The fact that $\cDes = \cDes'$ follows from Lemma~\ref{lem:phi-sufficient}.
\end{proof}

In fact, this argument shows that $(\phi, \cDes)$ and $(\phi', \cDes')$
are the same for a family of $\skewshape$ including both rectangles and strips.
\begin{cor}
    In any skew shape $\skewshape$ whose connected components are all rectangles,
    $(\phi, \cDes)$ and $(\phi', \cDes')$ are the same.
\end{cor}

\subsection{Hook Plus Internal Cell}

Let $\lambda = (n-k, 2, 1^{k-2})$ and $\mu = ()$.
The following construction is due to Adin, Elizalde, and Roichman.

\begin{defn}\cite[Corollary 4.8]{AER}\label{defn:hook-internal-cdes}
For $T\in \SYT(\skewshape)$,
let $n\in \cDes'(T)$ if and only if $\pos_T(T_{2,2}-1)$ is in the first column of $T$.
\end{defn}

This construction is different from the others in that it constructs $\psi = \phi'^{-1}$ instead of $\phi'$.
\begin{defn}\cite[Definition 4.13]{AER}\label{defn:hook-internal-phi}
    For $T\in \SYT(\skewshape)$, define $\psi T$ as follows.
    \begin{itemize}
        \item If $1\not\in \Des(T)$, let the first row of $\psi T$ have entries $[n]\setminus \cDes'(T)$,
        in increasing order from west to east, where $\cDes'(T)$ is defined in Definition~\ref{defn:hook-internal-cdes}.
        Let $(\psi T)_{2,2}$ be defined as follows, depending on the location of $\pos_T(T_{2,2}-1)$.
        \begin{enumerate}[label=(\alph*)]
            \item If $\pos_T(T_{2,2}-1)$ is in the first row of $T$, $(\psi T)_{2,2}=T_{2,2}-1$.
            \item If $\pos_T(T_{2,2}-1)$ is in the first column of $T$ and not the southernmost entry $(k,1)$,
            and $a$ is the entry in its southern neighbor, $(\psi T)_{2,2}=a-1$.
            \item If $\pos_T(T_{2,2}-1) = (k,1)$, $(\psi T)_{2,2}=n$.
        \end{enumerate}
        The remaining entries are in the first column of $T$, in increasing order from north to south.
        \item If $1\in \Des(T)$, let $\psi T = (\psi T^{\tra})^{\tra}$, reducing to the previous case.
    \end{itemize}
\end{defn}

\begin{thm}\label{thm:aer-hook-internal}\cite[Theorem 1.11]{AER}
    Let $\cDes'$ be defined by Definition~\ref{defn:hook-internal-cdes}, and let $\phi' = \psi^{-1}$,
    where $\psi$ is defined by Definition~\ref{defn:hook-internal-phi}.
    Then $(\cDes',\phi')$ is a cyclic descent map.
\end{thm}

Note this corollary of Lemma~\ref{lem:des-size-bounds}.
\begin{cor}\label{cor:cdes-unique}\cite[Lemma 4.1, Theorem 4.4]{AER}
    When $\skewshape = (n-k,2,1^{k-2})/()$,
    the minimum and maximum values of $|\Des(T)|$, for $T\in \SYT(\skewshape)$,
    are $k-1$ and $k$.
    Thus, in any cyclic descent map $(\cDes', \phi')$ of $\skewshape$,
    $|\cDes'(T)|$ has constant value $k$ over all $T\in \SYT(\skewshape)$,
    and is uniquely defined by letting $n\in \cDes'(T)$
    if and only if $|\Des(T)|=k-1$.
\end{cor}

\begin{lem}\label{lem:psi-first-row}
    Suppose $T\in \SYT(\skewshape)$ and $\Des(T)=S$, where
    $S\subseteq [n-1]$ and $|S|=k$.
    Then the first row of $T$ has entries $[n]\setminus (S+1)$.
\end{lem}
\begin{proof}
    If for any $x\in S$, $\pos_T(x+1)$ is in the first row of $T$,
    then $x$ cannot be a descent of $T$.
    Thus, all of $S+1$ is not in the first row of $T$.
    But $|S|=k$, so the $n-k$ elements of $[n]\setminus (S+1)$
    must all be in the first row of $T$.
\end{proof}

\begin{prop}
    The cyclic descent map $(\cDes', \phi')$ defined in Theorem~\ref{thm:aer-hook-internal}
    coincides with the $(\cDes, \phi)$ defined in Theorem~\ref{thm:main}.
\end{prop}
\begin{proof}
    Corollary~\ref{cor:cdes-unique} immediately implies that $\cDes = \cDes'$,
    because this function is unique.

    To show that $\phi = \phi'$, we will show that $\phi^{-1}$ coincides
    with the $\psi$ defined in Definition~\ref{defn:hook-internal-phi}.
    Since $\phi^{-1}$ commutes with transpose,
    and $\psi$ is defined for the case $1\in \Des(T)$ by
    $\psi T = (\psi T^{\tra})^{\tra}$,
    it suffices to handle the case $1\not\in \Des(T)$.

    Recall that
    \begin{equation*}
        \phi^{-1} = \RotSE^{-1} \circ \dem \circ \RotNW.
    \end{equation*}
    As $\RotNW$ is the identity map on $\SYT(\skewshape)$,
    \begin{equation*}
        \phi^{-1} = \RotSE^{-1} \circ \dem.
    \end{equation*}

    Suppose $1\not\in \Des(T)$.  Since $\cDes=\cDes'$,
    \begin{equation*}
        \cDes(\phi^{-1} T) = \cDes(\psi T) = \cDes(T)-1,
    \end{equation*}
    where this set does not contain $n$ because $1\not\in \Des(T)$.
    Denote this set $S$.
    Since $|\cDes(T)|=k$ by Corollary~\ref{cor:cdes-unique}, $|S| = k$ as well.
    By Lemma~\ref{lem:psi-first-row},
    the first row of both $\phi^{-1}(T)$ and $\psi T$ is
    \begin{equation*}
        [n]\setminus (S+1) = [n]\setminus \cDes(T).
    \end{equation*}
    Therefore, if we can show $(\psi T)_{2,2}=(\phi^{-1} T)_{2,2}$
    for all $T\in \SYT(\skewshape)$, we will be done.

    As $1\not\in \Des(T)$, $T_{1,2}=2$.  So, the demotion path of $T$ is either
    \begin{enumerate}[label=(\roman*)]
        \item a straight-east path from $(1,1)$ to $(1, n-k)$, or
        \item an L-shaped path from $(1,1)$ to $(1,2)$ to $(2,1)$.
    \end{enumerate}
    In the analysis below,
    we do casework on whether $\pos_T(T_{2,2}-1)$ is in case (a), (b), or (c)
    in Definition~\ref{defn:hook-internal-phi},
    and on which of (i) and (ii) is the demotion path of $T$.

    Case 1: $\pos_T(T_{2,2}-1)$ is in case (a).

    Note that $T_{2,2}-1\neq T_{1,2}$, because $T_{1,2}=2$ and $T_{2,2}\ge 4$.
    Thus, $\pos_T\lt(T_{2,2}-1\rt)$ is strictly due east of $(1,2)$.
    This implies that $T_{1,3}<T_{2,2}$, so the demotion path of $T$ is (i).
    So, $(\dem T)_{1,n-k} = n$.
    Note that demotion on $T$ decrements $T_{2,2}$ without moving it;
    moreover, it decrements $T_{2,2}-1$, which is in the first row, and moves it one cell westward.
    Thus, $\pos_{\dem T}(T_{2,2}-1) = (2,2)$ and $\pos_{\dem T}(T_{2,2}-2)$ is in the first row of $\dem T$.

    In $\dem T$, the set $\{x\in [n] | x\ge T_{2,2}-2\}$
    is southeast min-unimodal, with $T_{2,2}-1$ at $(2,2)$ and $T_{2,2}-2$ in the first row.
    By southeast min-unimodality,
    any additional elements in $\RcSE(\dem T)$ are also in the first row.
    Thus, the southeast rotation endpoint of $\dem T$ is $(1,n-k)$.
    Since $(\dem T)_{1, n-k} = n$, $\RotSE$ is the identity operator on $\dem T$.
    Therefore $\phi^{-1}T = \RotSE^{-1}(\dem T) = \dem T$,
    and $(\phi^{-1} T)_{2,2} = (\dem T)_{2,2} = T_{2,2}-1$, as desired.

    Case 2: $\pos_T(T_{2,2}-1)$ is in cases (b) or (c),
    and the demotion path of $T$ is (i).

    In this and the following case, we will determine
    $\phi^{-1} T = \RotSE^{-1}(\dem T)$
    by determining the action of $\RotSE^{-1}$ on $\dem T$.
    To do this we will determine the set $\RcSE(\phi^{-1} T)$ according
    to Lemma~\ref{lem:recover-rcse}.

    Since the demotion path is (i), $(\dem T)_{1,n-k} = n$.
    Since $(1,n-k)$ does not have a northern neighbor, in $\dem T$ the set
    $\RcSE(\phi^{-1} T)$ must be balanced at its western balance point, $(1,2)$.

    Since neither $T_{2,2}$ nor $T_{2,2}-1$ are in the first row of $T$,
    demotion on $T$ decrements these entries without moving them.
    Hence,
    \begin{equation*}
        \pos_{\dem T}(T_{2,2}-1) = (2,2),
        \quad
        \pos_{\dem T}(T_{2,2}-2) = \pos_T(T_{2,2}-1).
    \end{equation*}
    Note in particular that $\pos_{\dem T}(T_{2,2}-2)$ is not in the first row.
    In $\dem T$, the set
    \begin{equation*}
        S=\{x\in [n-1] | x\ge T_{2,2}-1\}
    \end{equation*}
    is southeast min-unimodal and balanced at $(1,2)$.
    Since $\pos_{\dem T}(T_{2,2}-2)$ is not in the first row,
    the set $S \cup \{T_{2,2}-2\}$ is not balanced at $(1,2)$ in $\dem T$.
    It follows that $\RcSE(\phi^{-1} T) = S$.

    Since $\pos_{\dem T}(T_{2,2}-2) = \pos_{\dem T}((\dem T)_{2,2}-1)$ is not in the first row, $(\dem T)_{2,2}-1 \neq (\dem T)_{1,2}$.
    So, $(\dem T)_{1,2} < (\dem T)_{2,2}-1 = T_{2,2}-2$, which implies $(1,2)\not\in \RpSE(\phi^{-1} T)$.
    So, the cells $\pos_{\dem T}(n)=(1,n-k)$ and
    \begin{equation*}
        \pos_{\dem T}(\min \RcSE(\phi^{-1} T))=\pos_{\dem T}(T_{2,2}-1)=(2,2)
    \end{equation*}
    are in different connected components of $\RpSE(\phi^{-1} T)$,
    and $\pos_{\dem T}(n)$ is northeast of $\pos_{\dem T}(\min \RcSE(\phi^{-1} T))$.
    Therefore, the action of $\RotSE^{-1}$ on $\dem T$ sends $n$ to the
    southwesternmost cell of $\RpSE(\phi^{-1} T)$,
    and the remaining elements of $\RcSE(\phi^{-1} T)$ northeast by one position.
    Note that all cells strictly due south of $\pos_T(T_{2,2}-1)=\pos_{\dem T}(T_{2,2}-2)$ have entries larger than $T_{2,2}-2$,
    and are thus in $\RpSE(\phi^{-1}(T))$.

    If $\pos_T(T_{2,2}-1)$ is not the southernmost cell in the first column,
    the entry rotated into $(2,2)$ is the entry immediately south of $(\dem T)_{2,2}-1$ in $\dem T$.
    Note that demotion on $T$ decrements both $T_{2,2}-1$ and the entry $a$ immediately to its south without moving them;
    thus, the entry rotated into $(2,2)$ by $\RotSE^{-1}$ on $\dem T$ is $a-1$, as desired.

    If $\pos_T(T_{2,2}-1)$ is the southernmost cell in the first column,
    there are no cells in $\RpSE(T)$ in the first column.
    Then, $(2,2)$ is the southwesternmost cell of $\RpSE(\phi^{-1}(T))$,
    and the entry rotated into $(2,2)$ is $n$, as desired.

    Case 3: $\pos_T(T_{2,2}-1)$ is in cases (b) or (c),
    and the demotion path of $T$ is (ii).

    Since the demotion path is (ii), $(\dem T)_{2,2}=n$ and
    $(\dem T)_{1,2} = T_{2,2}-1$.
    Moreover, because $T_{2,2} > T_{2,1}$, $(\dem T)_{1,2} > (\dem T)_{2,1}$.
    Thus, in $\dem T$, the set
    \begin{equation*}
        S=\{x\in [n-1] | x\ge T_{2,2}-1\}
    \end{equation*}
    is southeast min-unimodal and balanced at $(1,2)$.
    Moreover, because $T_{2,2}-1$ is in the first column of $T$,
    demotion of $T$ decrements this entry without moving it. So,
    \begin{equation*}
        \pos_T(T_{2,2}-1) = \pos_{\dem T}(T_{2,2}-2)
    \end{equation*}
    is in the first column of $T$.
    Thus, the set $S\cup \{T_{2,2} - 2\}$ is not balanced at $(1,2)$ in $\dem T$.
    It follows that $\RcSE(\phi^{-1} T) = S$.

    In this case, the cells $\pos_{\dem T}(n)=(2,2)$ and
    \begin{equation*}
        \pos_{\dem T}(\min \RcSE(\phi^{-1} T))=\pos_{\dem T}(T_{2,2}-1)=(1,2)
    \end{equation*}
    are in the same connected component of $\RpSE(\phi^{-1} T)$,
    and $\pos_{\dem T}(n)$ is southwest of $\pos_{\dem T}(\min \RcSE(\phi^{-1} T))$.
    So, $\RotSE^{-1}$ on $\dem T$ sends $n$ to the
    southwesternmost cell of $\RpSE(\phi^{-1} T)$,
    and the remaining elements of $\RcSE(\phi^{-1} T)$ northeast by one position.
    As in the previous case, all cells strictly due south of $\pos_T(T_{2,2}-1)=\pos_{\dem T}(T_{2,2}-2)$ are in $\RpSE(\phi^{-1}(T))$.

    If $\pos_T(T_{2,2}-1)$ is not the southernmost cell in the first column,
    the entry rotated into $(2,2)$ is the entry immediately south of $(\dem T)_{2,2}-1$ in $\dem T$.
    By the same reasoning as in the previous case, this is $a-1$, as desired.

    If $\pos_T(T_{2,2}-1)$ is the southernmost cell in the first column,
    there are no cells in $\RpSE(T)$ in the first column.
    As in the previous case, the entry rotated into $(2,2)$ is $n$, as desired.
\end{proof}

\subsection{Two-Row Young Diagrams}

Let $\lambda = (n-k,k)$ with $2\le k\le n/2$ and $\mu = ()$.
The following construction is due to Adin, Elizalde, and Roichman.

\begin{defn}\label{defn:two-row-cdes}\cite[Definition 1.12]{AER}
    For $T\in \SYT(\skewshape)$,
    let $n\in \cDes'(T)$ if and only if both of the following conditions hold.
    \begin{enumerate}
        \item $T_{2,k} = T_{2,k-1}+1$.
        \item For every $1<i<k$, $T_{2,i-1} > T_{1,i}$.
    \end{enumerate}
\end{defn}

By a slight abuse of notation, let $1+T_{\le x}$ denote the tableau obtained by
adding $1$ modulo $x$ to the entries $1,\ldots,x$ in $T$,
and leaving the remaining entries unchanged.

\begin{defn}\label{defn:two-row-phi}\cite[Definition 5.13]{AER}
    For $T\in \SYT(\skewshape)$, define $\phi' T$ as follows.
    \begin{enumerate}
        \item If $T_{2,k} = T_{2,k-1}+1$, let $\phi' T = \jdt(1+T_{\le x})$, where $x=T_{2,k}$.
        \item Otherwise:
        \begin{enumerate}
            \item If $n$ is in the first row of $T$, let $\phi' T = \pro(T)$.
            \item If $n$ is in the second row of $T$,
            let $\phi' T$ be the tableau obtained from $1+T$ as follows:
            switch $1$ with $y+1$, where $y$ is the entry immediately to the west of $1$, and then apply $\jdt$.
        \end{enumerate}
    \end{enumerate}
\end{defn}

\begin{thm}\label{thm:aer-two-row}\cite[Theorem 1.14]{AER}
    The $(\cDes',\phi')$ defined by Definitions~\ref{defn:two-row-cdes}
    and \ref{defn:two-row-phi} is a cyclic descent map.
\end{thm}

\begin{lem}\label{lem:two-row-se-rot-endpt}
    Let $T\in \SYT(\skewshape)$.
    The southeast rotation endpoint of $T$ is $(2,k)$ if and only if
    one of the following conditions hold.
    \begin{itemize}
        \item $T_{2,k}=T_{2,k-1}+1$.
        \item $T_{2,k}=n$ and $T_{1,k} = T_{2,k-1} + 1$.
    \end{itemize}
\end{lem}
\begin{proof}
    Suppose $T_{2,k}=T_{2,k-1}+1$, and let $T_{2,k}=x$.
    Then the easternmost entries in the first row of $T$ are
    $x+1,x+2,\dots,n$ (where this list is empty if $n=x$),
    and the easternmost two entries in the second row of $T$ are $x,x-1$.
    By southeast min-unimodality of $\RcSE(T)$,
    $\min \RcSE(T)$ is in the second row of $T$,
    so the southeast rotation endpoint of $T$ is $(2,k)$.

    Else, suppose $T_{2,k}\neq T_{2,k-1}+1$.
    We consider two cases: $T_{2,k}=n$ and $T_{2,k}\neq n$.

    If $T_{2,k}=n$, let $T_{2,k-1}=z<n-1$,
    so the easternmost entries in the first row of $T$ are $z+1,\dots,n-1$.
    Thus, $\RcSE(T)=\{z+1,\dots,n\}$.

    Since $\pos_T(z) = (2,k-1)$ and $z+1>z$, $\pos_T(z+1)$ is nonstrictly due east of $(1,k)$.
    We do casework on whether $\pos_T(z+1) = (1,k)$.
    If $\pos_T(z+1) = (1,k)$, then the southeast rotation endpoint of $T$ is
    $(2,k)$.  Moreover
    \begin{equation*}
        T_{1,k} = z+1 = T_{2,k-1}+1 ,
    \end{equation*}
    as desired.
    Otherwise, $\pos_T(z+1)$ is strictly east of $(1,k)$.
    Then the southeast rotation endpoint of $T$ is $(1,n-k)$.
    Since $\pos_T(z+1)$ is strictly east of $(1,k)$,
    \begin{equation*}
        T_{1,k} < z+1 = T_{2,k-1},
    \end{equation*}
    so neither of the two listed conditions holds, as desired.

    If $T_{2,k}\neq n$, let $T_{2,k}=x$.
    The easternmost entries in the first row of $T$ are $x-1,x+1,x+2,\dots,n$,
    while the easternmost entry in the second row of $T$ is $x$.
    By southeast min-unimodality, $\min(\RcSE(T))$ is in the first row of $T$,
    and the southeast rotation endpoint of $T$ is $(1,n-k)$.
    Note that $\pos_T(\min \RcSE(T))$ may be the southeast interior corner $(1,k)$,
    but because $\pos_T(n)=(1,n-k)$ is northeast of $(1,k)$,
    the southeast rotation endpoint of $T$ is still $(1,n-k)$.
    As neither of the two listed conditions holds, the claim is proved.
\end{proof}

Note that $\RotNW$ is the identity map on $\SYT(\skewshape)$, so
\begin{equation*}
    \phi = \pro \circ \RotSE.
\end{equation*}

\begin{prop}
    The cyclic descent map $(\cDes', \phi')$ defined in Theorem~\ref{thm:aer-two-row}
    coincides with the $(\cDes, \phi)$ defined in Theorem~\ref{thm:main}.
\end{prop}
\begin{proof}
    Let $T_{2,k}=x$. We do casework on whether $T_{2,k} = T_{2,k-1}+1$.

    Case 1: $T_{2,k} = T_{2,k-1}+1$.

    By Lemma~\ref{lem:two-row-se-rot-endpt},
    the southeast rotation endpoint of $T$ is $(2,k)$.
    We do casework on whether $x=n$.

    If $x=n$, then $T_{2,k}=n$ and $T_{2,k-1}=n-1$.
    By southeast min-unimodality, $\RcSE(T)$ is contained in the second row of $T$.
    Thus, the southeast rotation endpoint of $T$ is $(2,k)=\pos_T(n)$.
    So, the action of $\RotSE$ on $T$ is the identity, and $\phi T = \pro T$.
    As
    \begin{equation*}
        \phi' T = \jdt(1+T_{\le n}) = \pro T,
    \end{equation*}
    we have $\phi T = \phi' T$.

    If $x\neq n$, then the easternmost entries in the first row of $T$ are $x+1,x+2,\dots,n$,
    and the easternmost entries in the second row are $x-1, x$.
    The set $\{x-1,\ldots,n\}$ is southeast min-unimodal in $T$,
    and thus contained in $\RcSE(T)$; moreover,
    by southeast min-unimodality, all entries in $\RcSE(T)$ less than $x-1$ are in the second row.
    Thus the southeast rotation endpoint of $T$ is $(2,k)$,
    and $\RotSE$ rotates the entries $\{x,\ldots,n\}$.

    In $\RotSE(T)$, the easternmost entries in the first row are $x,x+1,\dots,n-1$,
    and $n$ is in $(2,k)$.
    Consequently, the tableau $1+\RotSE(T)$ is identical to $1+T_{\le x}$,
    and
    \begin{equation*}
        \phi T = \pro(\RotSE(T)) = \jdt(1+\RotSE(T)) = \jdt(1+T_{\le x}) = \phi' T,
    \end{equation*}
    as desired.

    Case 2: $T_{2,k}\neq T_{2,k-1}+1$.

    We again do casework on whether $x=n$.

    If $x=n$, let $T_{2,k-1} = z < n-1$.
    The easternmost entries in the first row of $T$ are $z+1,\dots,n-1$.

    If $T_{1,k}=z+1$, then by Lemma~\ref{lem:two-row-se-rot-endpt},
    $\RotSE$ does not alter $T$.  Thus, $\phi T = \pro T$.
    $\phi' T$ is obtained from $1+T$ by swapping $1$ with $z+2$,
    and then applying $\jdt$.
    Since $\pos_{1+T}(z+2)=(1,k)$ is the northern neighbor of $\pos_{1+T}(1)=(2,k)$,
    this is equivalent to just promoting $T$.
    So, $\phi' T = \pro T = \phi T$.
    Otherwise, if $T_{1,k} \neq z+1$, then $\RotSE$ rotates the set $\{z+1,\dots,n\}$.
    In $\RotSE(T)$, the easternmost entries in the first row are $z+2,\dots,n$,
    and $z+1$ is in $(2,k)$.
    Promoting this yields the same result as swapping $1$ with $z+2$ in $1+T$
    and applying $\jdt$.  Hence, $\phi T = \phi' T$.

    If $x\neq n$, the easternmost entries in the first row of $T$ are $x-1,x+1,x+2,\ldots,n$.
    The set $\{x-1,\ldots,n\}$ is southeast min-unimodal in $T$,
    and thus contained in $\RcSE(T)$; moreover,
    by southeast min-unimodality, all entries in $\RcSE(T)$
    smaller than $x-1$ are in the first row.
    Thus the southeast rotation endpoint of $T$ is $(1,n-k)$,
    and $\RotSE$ does not alter $T$.
    So $\phi T = \pro T$.
    As $\phi' T = \pro T$ by definition, we have $\phi T = \phi' T$.

    Therefore, $\phi$ and $\phi'$ coincide in all cases.
    By Lemma~\ref{lem:phi-sufficient}, $\cDes$ and $\cDes'$ coincide as well.
\end{proof}

\subsection{Young Diagram with Disconnected Northeast Cell}

The following construction is due to Elizalde and Roichman.

\begin{thm}\cite[Proposition 5.3]{ER}\label{thm:er}
    Let $\skewshape$ consist of a Young diagram of size $n-1$ and a disconnected northeast cell.
    For $T\in \SYT(\skewshape)$, let $n\in \cDes(T)$ if and only if one of the following holds.
    \begin{enumerate}
        \item $n$ is in the disconnected cell.
        \item $d\neq n$ is in the disconnected cell and $n-d\in \Des(\jdt(-d+T))$
    \end{enumerate}
    Moreover, define $\phi'$ by
    \begin{equation*}
        \phi'(T) = \jdt(1+d+\jdt(-d+T)).
    \end{equation*}
    Then $(\cDes',\phi')$ is a cyclic descent map.
\end{thm}

This cyclic descent map is \textit{not} a special case of $(\cDes,\phi)$.
As a counterexample, let
\begin{equation*}
    T = \young(:::1,235,46).
\end{equation*}
Then
\begin{equation*}
    \phi T = \young(:::2,134,56)
    \quad
    \text{and}
    \quad
    \phi' T = \young(:::2,146,35),
\end{equation*}
so $\phi T \neq \phi' T$.
Moreover, $\cDes(\phi T) = \{2,4,6\}$, while $\cDes'(\phi T) = \{2,4\}$.

\section{Concluding Remarks and Open Problems}\label{sec:conclusion}

With the results in this paper, we have fully solved the problem of
finding an explicit cyclic descent map for all skew shapes where a cyclic descent map exists.
We thus have a constructive proof of Theorem~\ref{thm:arr}.

As noted in Remark~\ref{rem:not-zn-action},
the $\phi$ defined in Theorem~\ref{thm:main} does not, in general, generate a $\ZZ_n$-action.
As promotion generates a $\ZZ_n$ action on rectangles,
this is one property of promotion on rectangles that does not generalize to $\phi$.
It would be interesting to determine when this property holds.

\begin{prob}
    For which $\skewshape$, not a connected ribbon, does $\phi$ generate a $\ZZ_n$-action?
\end{prob}
This class includes rectangles \cite[Theorem 4.1(a)]{StaPE} and strips \cite[Proposition 3.3]{AER};
it isn't difficult to show that it includes any skew shape
whose connected components are all rectangles.
\cite[Theorem 1.11]{AER} shows that hooks plus an interior cell are in this class as well.

It would also be interesting to study the structure of orbits of the action of $\phi$.
As $\phi$ does not generate a $\ZZ_n$-action, these orbit sizes are not always divisors of $n$.
In fact, it is possible for an orbit size to be neither a multiple nor a divisor of $n$;
as noted in \cite{AER}, the orbit generated by $\phi$ on
\begin{equation*}
    T=\young(13479,2568)
\end{equation*}
has size $6$ and the period of the corresponding cyclic descent sets is $3$.

Moreover, computer experiments done by the author show that the orbit generated by $\phi$ on
\begin{equation*}
    T = \young(1259\fourteen,378\ten\nineteen,4\twelve\fifteen\eighteen\twenty,6\seventeen\twthree\twfour\twsix,\eleven\twone\twfive,\thirteen\twtwo\twseven,\sixteen\tweight\twnine)
\end{equation*}
has size $488969=29\times 16861$,
suggesting that these orbit sizes can grow arbitrarily large, and rather quickly.

The following problems aim to better understand the distribution of orbit sizes.

\begin{prob}
    For $\skewshape$ not a connected ribbon, determine, exactly or asymptotically
    in terms of $n$ and the row lengths of $\lambda$ and $\mu$,
    the number of distinct orbits of $\phi$ on $\SYT(\skewshape)$.
\end{prob}

\begin{prob}
    For $\skewshape$ not a connected ribbon, determine, exactly or asymptotically
    in terms of $n$ and the row lengths of $\lambda$ and $\mu$,
    the size of the largest orbit of $\phi$ on $\SYT(\skewshape)$.
\end{prob}

\end{document}